\documentclass{article}

\usepackage[english]{babel}

\usepackage[papersize={8.5in,11in}, margin=1in,footskip=0.4in,vmarginratio={1:1},hmarginratio={1:1}]{geometry}

\usepackage{amsfonts} 
\usepackage{amsmath} 
\usepackage{amssymb} 
\usepackage{amsthm} 
\usepackage{latexsym} 
\usepackage{mathrsfs} 
\usepackage{mathtools}
\usepackage{xcolor}
\usepackage{enumitem}

\numberwithin{equation}{section}

\newtheorem{theorem}{Theorem}[section]
\newtheorem{remark}{Remark}[section]
\newtheorem{lemma}{Lemma}[section]
\newtheorem{prop}{Proposition}[section]

\newcommand{\xbar}[0]{\overline{x}}
\newcommand{\parbar}[0]{\overline{\partial}}
\newcommand{\japbr}[0]{\langle\overline \partial \rangle}

\newcounter{1}
\newcounter{2}
\newcounter{3}
\newcounter{4}
\title{Local Well-posedness of the Free-boundary Problem in Incompressible Elastodynamics with Surface Tension}
\date{}
\author{Longhui Xu}

\begin{document}
	\maketitle
	\begin{abstract}
	We prove the local well-posedness of the 3D free-boundary incompressible elastodynamics with surface tension describing the motion of an elastic medium in a periodic domain with a moving graphical surface. The deformation tensor is assumed to satisfy the neo-Hookean linear elasticity. We adapt the idea in \cite{zhang2024well} to generate an approximate problem with artificial viscosity indexed by $\kappa>0$ to boost the boundary regularity, which recovers the original system as $\kappa \to 0,$ and the energy estimates yield no regularity loss. 
	\end{abstract}
\tableofcontents
\pagebreak

\section{Introduction}
Let $x = (x_1,x_2,x_3)$ and $\xbar = (x_1,x_2)$. We define $\mathcal{D}_t:=\{\xbar\in \mathbb{T}^2, -b<x_3<\psi(t,\xbar)\}$ to be the periodic domain at time $t$ with boundary described by a graph $\psi$ and $ \mathcal{D}:= \bigcup _{0\leq t\leq T} \{t\} \times \mathcal{D}_t $. We consider the incompressible free-boundary elastodynamic equations
	
	\begin{equation} \label{inisys}
		\begin{cases}
			D_tu +\nabla p=\text{div}(\mathbf{F}\mathbf{F}^T)  \qquad & \text{in } \mathcal{D},\\
			\text{div }u=0 \qquad & \text{in } \mathcal{D},\\
			D_t \mathbf{F}=\nabla u \mathbf{F} \qquad & \text{in } \mathcal{D},\\
			\text{div }\mathbf{F}^T=0 \qquad & \text{in } \mathcal{D}.
		\end{cases}
	\end{equation}
	 The velocity and fluid pressure of the elastic medium are denoted by $u,p$. $\mathbf{F}:=(\mathbf{F}_{ij})_{3\times3}$ is the deformation tensor and $\mathbf{F}^T:=(\mathbf{F}_{ji})$ is the transpose of $F$. $\mathbf{F}\mathbf{F}^T$ is the Cauchy-Green tensor for neo-Hooken elsatic materials. $\nabla := (\partial_1,\partial_2,\partial_3)$ is the standard spatial derivative, div $X:=\nabla \cdot X= \partial_i X^i$ is the standard divergence for any vector filed $X$. $D_t:= \partial_t +v\cdot \nabla$ is the material derivative.
	 
	 We consider the boundary conditions
	 \begin{equation}\label{inibd}
		\begin{cases}
			D_t|_{\partial \mathcal {D}}\in\mathcal{T}(\partial \mathcal{D}), \\
			p=\sigma \mathcal {H}  \qquad & \text{on } \partial\mathcal{D}_t,\\
			\mathbf{F}^T N =0 \qquad & \text{on } \partial \mathcal{D}_t.
		\end{cases}
	\end{equation}
	where $\mathcal{T}(\partial \mathcal {D}_t)$ is the tangent bundle of $\partial \mathcal {D}_t$, $\sigma$ is the constant coefficient of surface tension, $\mathcal {H}$ is the mean curvature of the surface and $N$ is the normal vector to the boundary. The first condition indicates that the boundary moves with the velocity of the fluid, and the second condition shows that the pressure is balanced by surface tension. The third condition originates from the Rankine-Hugoniot conditions for the vortex sheets in compressible elastodynamics, which leads to $\det F =0$ on the boundary \cite{Zhang_2022,chen2017linear,chen2020linear,chen2020nonlinear,chen2023stabilization}. 
	\begin{remark}[Initial constraints for the deformation tensor]
		The last conditions of (\ref{inisys}) and (\ref{inibd}) are only required to hold for the initial data and they propagate within the lifespan of the solutions. So the system is not overdetermined with these two conditions. 
	\end{remark}
	
	Let $\mathbf F_j$ be the $j$-th column of $\mathbf F$. Given the initial data $(u_0, \mathbf F_k^0, p_0, \psi_0)$, we study the local well-posedness of the following system for the case $\sigma>0$
	\begin{equation}\label{sys}
		\begin{cases}
			D_tu +\nabla p= (\mathbf{F}_k\cdot \nabla) \mathbf{F}_k \qquad & \text{in } \mathcal{D},\\
			\text{div }u=0 \qquad & \text{in } \mathcal{D},\\
			D_t\mathbf{F}_j=( \mathbf F_j\cdot \nabla) u \qquad & \text{in } \mathcal{D},\\
			\text{div }\mathbf{F}_j=0 \qquad & \text{in } \mathcal{D},\\
			D_t|_{\partial \mathcal {D}}\in\mathcal{T}(\partial \mathcal{D}), \\
			\mathbf{F}_j\cdot N =0 \qquad & \text{on } \partial \mathcal{D}_t, \\
			p=\sigma \mathcal {H}  \qquad & \text{on } \partial\mathcal{D}_t,\\
			(u,\mathbf F_k, p, \psi)|_{t=0} = (u_0, \mathbf F_k^0, p_0, \psi_0).
		\end{cases} 
	\end{equation}

\subsection{Reformulation in flattened domains}
One way to convert (\ref{inisys}) to a system of equations on a fixed domain is by Lagrangian coordinates, it is powerful in the sense that it can be used to solve problems with a general domain, instead of the domain we considered here with boundary described by a graph. Also, we do not have to worry about the propagation of the initial constraints of the deformation tensor $\mathbf F$ when constructing the approximate system, since the equations reformulated in Lagrangian coordinates automatically only depend on the initial data of  $\mathbf F$. However, it leads to painful calculations when dealing with surface tension. Here we introduce graphical coordinates to convert the free-boundary problem into a fixed-domain problem on 
	\begin{equation*}
		\Omega:=\mathbb{T}^2 \times (-b,0).
	\end{equation*}
	We define a family of diffeomorphisms $\displaystyle \Phi(t,\cdot) :\Omega \to \mathcal{D}_t$ by 
	\begin{align}
	\Phi(t,\overline{x},x_3) = (\overline{x},\varphi(t,\overline{x},x_3)),
	\end{align}
	where $\varphi(t,\cdot)$ is an extension of the boundary graph $\psi$ to the interior by 
	\begin{align}
	\varphi(t,\overline{x},x_3)=x_3+\chi(x_3)\psi(t,\overline{x})	
	\end{align}
	and $\chi\in C_c^\infty(-b,0]$ satisfying 
	\begin{equation}\label{chi}
		||\chi'||_{L^\infty (-b,0]}\leq \frac{1}{b+|\psi_0|_\infty}, \qquad \chi\equiv 1 \text{ on }(-\delta_0,0]
	\end{equation}
	for some small constant $\delta_0>0$. To avoid the case that the moving surface touches the bottom, we can assume $|\psi_0|_\infty< \frac{b}{3}$, it can be shown that $\sup_{[0, T_0]}|\psi(t)|_\infty<\frac{b}{2}$ for some $T_0>0$ with 
	\begin{align}
	\partial_3 \varphi (t,\overline{x},x_3) =1+\chi'(x_3)\psi(t,\overline{x})\geq 1-\frac{|\psi(t)|_\infty}{b+|\psi_0|_\infty}\geq 1-\frac{\frac{b}{2}}{b+\frac{b}{3}} =: c_0>0,  \qquad t\in[0,T], \label{partial3phi}	
	\end{align}
	which ensures $\Phi(t)$ to be well-defined diffeomorphisms. 
	
	Now we introduce the converted variables
	\begin{equation}
		v(t,x)=u(t,\Phi(t,x)), \qquad q(t,x)=p(t,\Phi(t,x)), \qquad F(t,x)=\mathbf{F}(t,\Phi(t,x))m,
	\end{equation}
	the induced differential operators
	\begin{align}
		\partial_t^\varphi&=\partial_t-\frac{\partial_t \varphi}{\partial_3 \varphi}\partial_3,\\
		\nabla_\tau^\varphi=\partial_\tau^\varphi&=\partial_\tau-\frac{\partial_\tau \varphi}{\partial_3 \varphi}\partial_3, \qquad \tau=1,2,\\
		\nabla_3^\varphi= \partial_3^\varphi&=\frac{1}{\partial_3 \varphi}\partial_3
	\end{align}
	and the material derivative 
	\begin{align}
		D_t^\varphi&:=\partial_t^\varphi+v\cdot\nabla^\varphi \\
				&=\partial_t+\overline{v}\cdot\overline{\partial}+\frac{1}{\partial_3 \varphi}(v\cdot\mathbf{N} -\partial_t\varphi)\partial_3 ,
	\end{align}
	where $\overline{\partial}=\overline{\nabla}=(\partial_1,\partial_2)$ and $\mathbf{N}:=(-\partial_1\varphi,-\partial_2\varphi,1)$ is the extension of the normal vector $N$ to the interior. 
	
	Let $\Sigma = \mathbb T^2\times \{0\}$ and $\Sigma_b = \mathbb T^2\times \{-b\}$. The system (\ref{sys}) can be reformulated as follows, 
		\begin{equation}\label{phisys}
		\begin{cases}
			D_t^\varphi v +\nabla^\varphi q= (F_k\cdot \nabla^\varphi) F_k \qquad & \text{in } \Omega,\\
			\nabla^\varphi \cdot v=0 & \text{in } \Omega,\\
			D_t^\varphi F_j=(F_j\cdot \nabla^\varphi) v  & \text{in } \Omega,\\
			\nabla^\varphi \cdot  F_j=0 \qquad & \text{in } \Omega,\\
			\partial_t \psi=v\cdot N &\text{on } \Sigma,\\
			q=-\sigma \overline{\nabla}\cdot\big(\frac{\overline{\nabla\psi}}{\sqrt{1+|\overline{\nabla}\psi|^2}}) &\text{on } \Sigma,\\
			F_j\cdot N=0 &\text{on } \Sigma,\\
			v_3=F_{3j}=0 &\text{on } \Sigma_b,\\
			(v,F_k, q, \psi)|_{t=0} = (v_0, F_k^0, q_0, \psi_0).
		\end{cases}
	\end{equation}

\subsection{History and background}
	The system of incompressible elastodynamics is a coupling system of the Euler equations with a transport equation for the deformation tensor. Let us first briefly review the results of the free-boundary incompressible Euler equations. The first breakthrough was by Wu \cite{wu1997well,wu1999well} on the local well-posedness of the irrotational water wave. For the general incompressible problem with non-zero vorticity, Christodoulou-Lindblad \cite{christodoulou2000motion} established an apriori energy estimate and Lindblad \cite{lindblad2001well,lindblad2005well} proved the LWP for the case without surface tension by Nash-Moser iteration. Coutand-Shkoller \cite{coutand2007well,coutand2010simple} introduced tangential smoothing to prove the LWP for the case with surface tension and avoided the loss of regularity. We also refer to \cite{alazard2014cauchy,zhang2008free} for the case without surface tension and \cite{schweizer2005three,shatah2008geometry,shatah2008priori,shatah2011local} for the case with nonzero surface tension.
	
	The study of elastodynamics has a long history. The global existence of the fixed-boundary problem is expected and we refer to  \cite{lei2005global,lei2008global,lin2005hydrodynamics,lin2008initial,lin2012some, coutand2006interaction,dafermos2005hyperbolic}. However, the free-boundary problem is more difficult as we have limited boundary regularity and the boundary terms can enter the highest order in the energy estimate. Most of the existing literature neglected the effect of surface tension and the fluid tends to be unstable. Gu-Wang \cite{gu2020well} and Li-Wang-Zhang \cite{li2021well} proved the LWP under a mixed stability condition. Hu-Huang \cite{hu2019well} proved the LWP under the Rayleigh-Taylor sign condition. Moreover, the free boundary problem can be considered as the one-phase problem of the vortex sheets in incompressible elastodynamics, we refer to Li-Wang-Zhang \cite{li2019well}. For the case with surface tension, Gu-Lei \cite{gu2020local} proved the local well-posedness under the force balance condition by studying the viscoelastic system in Lagrangian coordinates and the vanishing viscosity limit. 
	
		It is painful to study the case with nonzero surface tension with Lagrangian coordinates due to its complicated boundary description. In the presenting manuscript, we prove the local well-posedness of the free-boundary incompressible elastodynamics with surface tension under the boundary condition $(\ref{inibd})_3$ in graphical coordinates and the energy estimate without regularity loss.

\subsection{Main result}
The theorem concerns the local well-posedness of the free boundary incompressible elastodynamics modelled by (\ref{phisys}). The initial data $(\psi_0,q_0,v_0,F_k^0)$ satisfies the $m$-th order compatibility condition if, for $0\leq j\leq m$,
\begin{align}
\partial_t ^j q\big|_{t=0} =\sigma \partial_t^j\mathcal{H}\big|_{t=0} \quad&\text{on } \Sigma,\nonumber\\
  \partial_t^{j+1}\psi\big|_{t=0}=\partial_t^j	(v\cdot N)\big|_{t=0} \qquad &\text{on } \Sigma,\\
\partial_t^jv_3\big|_{t=0}=0 \quad &\text{on } \Sigma_b\nonumber.
\end{align}
\begin{theorem}
	Fix $\sigma>0$. Let $(q_0,v_0,F_k^0)\in  H^4(\Omega) $ and $\psi_0\in H^{5.5}(\Sigma)$ be initial data of the system (\ref{phisys}) that verifies the compatibility condition up to $3$-th order, $E(0)\leq M$ for some constant $M>0$, and the initial constraints 
	\begin{align}
		\nabla^{\varphi_0} \cdot  F^0_j=0 &\text{ in }\Omega,\\
		F^0_j\cdot N_0=0 &\text{ on } \Sigma\cup\Sigma_b.
	\end{align}
	Then there exists $T>0$ depending only on the initial data and $\sigma$ such that (\ref{phisys}) admits a unique solution $(v(t),F_k(t),\psi(t))$ satisfying the energy estimate 
	\begin{align}
	\sup_{[0\leq t\leq T]} E(t)\leq C(\sigma^{-1}) P(E(0)),	
	\end{align}
	where $P(\cdot )$ is a generic polynomial in its arguments and the energy is defined by 
	\begin{align}
		E(t)= \sum_{l=0}^4\bigg(\sum_{k=0}^3||\partial_t^l(v,F_k)||_{4-l}+|\sqrt\sigma\partial_t^l\psi|_{5-l}\bigg).
	\end{align}
	Here $||\cdot||_s$ represents represents the interior Sobolev norm $||\cdot ||_{H^s(\Omega)}$ and $|\cdot|_s$ represents the boundary Sobolev norm $|\cdot|_{H^s(\Sigma)}$.

\end{theorem}

\section{Uniform Estimates of the nonlinear approximate system}
The local well-posedness of the free-boundary fluid is not a direct consequence of the priori estimate \cite{luo2022compressible, Zhang_2022}. There is a loss of boundary regularity when trying to obtain the higher-order energy estimate of the linearized system for Picard iteration since we no longer have a delicate symmetry when dealing with the surface tension as in the priori energy estimate \cite{xu2024priori}.
One way is to use tangential smoothing to boost the regularity of $\psi$ \cite{luo2022compressible}, but the initial constraint $F_j\cdot N =0 $ would no longer propagate due to the mismatched terms arising from the tangentially smoothed $\psi$ when deriving the linear equation of $F\cdot N$ on the boundary. Here we use a simpler approximate scheme to boost the regularity of $\psi$ and $\psi_t$ introduced in \cite{zhang2024well} by only adding artificial viscosity terms to the balanced pressure boundary condition: 
\begin{align}
q=-\sigma \overline{\nabla}\cdot\big(\frac{\overline{\nabla\psi}}{\sqrt{1+|\overline{\nabla}\psi|^2}}) + \kappa(1-\overline{\Delta})^2\psi + \kappa(1-\overline{\Delta})\partial_t\psi.
\end{align}
This approximate scheme does not affect the propagation of the initial boundary constraint $F\cdot N = 0$, which only depends on the equation $D_t^\varphi F_j=(F_j\cdot \nabla^\varphi) v$ and the kinematic boundary condition $\partial_t \psi=v\cdot N$. 
Now we introduce the approximate system of (\ref{sys}) indexed by $\kappa > 0$,
		\begin{equation}\label{approxsys}
		\begin{cases}
			D_t^\varphi v +\nabla^\varphi q= (F_k\cdot \nabla^\varphi) F_k \qquad & \text{in } \Omega,\\
			\nabla^\varphi \cdot v=0 & \text{in } \Omega,\\
			D_t^\varphi F_j=(F_j\cdot \nabla^\varphi) v  & \text{in } \Omega,\\
			\nabla^\varphi \cdot  F_j=0 \qquad & \text{in } \Omega,\\
			\partial_t \psi=v\cdot N &\text{on } \Sigma,\\
			q=-\sigma \overline{\nabla}\cdot\big(\frac{\overline{\nabla\psi}}{\sqrt{1+|\overline{\nabla}\psi|^2}}) + \kappa(1-\overline{\Delta})^2\psi + \kappa(1-\overline{\Delta})\partial_t\psi &\text{on } \Sigma,\\
			F_j\cdot N=0 &\text{on } \Sigma,\\
			v_3=F_{3j}=0 &\text{on } \Sigma_b,\\
			(v,F_k, q, \psi)|_{t=0} = (v_0, F_k^0, q_0, \psi_0).
		\end{cases}
	\end{equation}
The energy functional associated with system \ref{approxsys} is defined by
\begin{align}
E_4^\kappa(t) = \sum_{l=0}^4||\partial_t^l( v, F_k)||^2_{4-l}+\sum_{l=0}^4(|\sqrt\sigma\partial_t^l \psi|_{5-l},|\sqrt\kappa \partial_t^l \psi|_{6-l}+\int_0^t|\sqrt\kappa \partial_t^{l+1} \psi|_{5-k}).
\end{align}
\begin{remark}[Reduction of pressure] The control of the pressure $q$ can be reduced to the control of $u$ and $F_k$. Invoking the momentum equation, for $\tau =1 ,2,$ we have 
	\begin{align}
		\partial_3 q &=\partial_3\varphi\partial_3^\varphi q= -\partial_3 \varphi D_t^\varphi v +\partial_3 \varphi(F_k\cdot \nabla^\varphi )F_{3k},\\
			\partial_\tau q	&=\partial_\tau^\varphi q+\frac{1}{\partial_3 \varphi}\partial_3 q=D_t^\varphi v_\tau -(F_k\cdot \nabla^\varphi )F_{\tau k} +\frac{1}{\partial_3 \varphi}\partial_3 q.
	\end{align}
	It follows that, for $i=1,2,3,s=1,2,3,$
	\begin{align*}
	||\partial_t^s	 \partial_i q||_{3-s}\leq P\big(\sum_{l=0}^3 |\partial_t ^l \overline{\nabla}\psi|_{3-l},\sum _{l=0}^4 ||\partial_t^l v||_{4-l}, \sum _{l=0}^4 ||\partial_t^k F||_{4-l}\big).
	\end{align*}
	We can obtain the $L^2$ estimate with Poincar\'e's inequality
	\begin{equation}
		||q||_0
		\lesssim _{vol(\Omega)} ||\partial q||_0 +\int_\Omega  q,
	\end{equation}
	where 
	\begin{equation}
		\int_\Omega q dx= \int_\Omega \partial_3 x_3 q dx= -\int_\Omega x_3\partial_3 q dx\leq ||x_1||_0||\partial q||_0\lesssim_{Vol(\Omega)}||\partial q||_0.
	\end{equation}
	The boundary terms vanishes since $x_3 = 0 $ on $\Sigma$ and $q=0$ on $\Sigma_b.$

\end{remark}
We aim to establish the a priori estimates of (\ref{approxsys}) uniform in $\kappa>0$ so that we can recover a solution to the original system (\ref{phisys}) by taking the limit $\kappa \to 0^+$.
\begin{prop} There exists some $T_\sigma >0$ independent of $\kappa$ such that, for $t\in[0,T_\sigma]$,
\begin{align}
E^\kappa_4(t)\leq C(\sigma^{-1})P(E_4^\kappa(0)).
\end{align}
\end{prop}


Let us introduce the following Reynold's transport theorems that are used throughout the analysis without further mentioning. 
\begin{lemma}[Integration by parts]\label{ibp}
		 Let $g=g(t,x), f=f(t,x), x\in \Omega$.
		\begin{align}
		\int_\Omega 	(\partial_i^\varphi f)g\partial_3 \varphi =-\int_\Omega 	 f(\partial_i^\varphi g)\partial_3 \varphi +\int_\Sigma fgN_i +\int_{\Sigma_b} fgn_i.
		\end{align}
	where $N=(-\partial_1\psi,-\partial_2 \psi ,1) ,n=(0,0,1)$ are the normal vector function of $\Sigma,\Sigma_b$ respectively.
	\end{lemma}
	
	\begin{lemma}[Transport Theorem]\label{transp}
		 Let $f=f(t,x), x\in \Omega$.
		\begin{align}
		\frac{d}{dt}\int_\Omega f\partial_3 \varphi =\int_\Omega D_t^\varphi f\partial_3 \varphi .
		\end{align}
	\end{lemma}
	The proofs are straightforward applications of change of variables and integration by parts, we refer to \cite{xu2024priori} for details.
\subsection{$L^2$ energy conservation}\label{L2}
\begin{prop}
	Let
	\begin{align}
	E_0^\kappa = \frac{1}{2}\int_\Omega\bigg(|v|^2 + \sum_{k=1}^3\int_\Omega |F_k|^2\bigg)\partial_3 \varphi	+\frac{1}{2}\int_\Sigma \bigg(\sigma\sqrt{1+|\overline{\nabla}\psi|^2}+\kappa|(1-\overline{\Delta})\psi |^2\bigg)+\int_0^t\int_\Sigma \kappa |\japbr\partial_t\psi|^2.
	\end{align}
	Here $\langle\cdot\rangle$ is the Japanese bracket, then 
	\begin{align}
	\frac{d}{dt}E_0^\kappa(t)= 0	.
	\end{align}
	within the lifespan of the solution to (\ref{approxsys}).
\end{prop}

\noindent\textit{Proof.} Testing the first equation of (\ref{approxsys}) with  $v\partial_3\varphi$, by integration by parts (\ref{ibp}) and transport theorem (\ref{transp}), the left-hand side gives
\begin{align*}
	\int_\Omega (D_t^\varphi v +\nabla^\varphi q) \cdot v\partial_3 \varphi = \frac{1}{2}\frac{d}{dt}\int_\Omega |v|^2\partial_3 \varphi +  \int_\Sigma q (v\cdot N),
\end{align*}
where the integral on $\Sigma_b$ vanishes due to the slip boundary conditions.
Integrating $(F_k\cdot \nabla^\varphi)$ by parts and invoking the third equation of (\ref{approxsys}), the right-hand side gives
\begin{align*}
	\int_\Omega (F_k\cdot \nabla^\varphi) F_k \cdot v\partial_3 \varphi 
		&= -\int_\Omega  F_k \cdot (F_k\cdot \nabla^\varphi) v\partial_3 \varphi
		=-\int_\Omega  F_k \cdot D_t^\varphi F_k \partial_3 \varphi =-\frac{1}{2}\frac{d}{dt}\sum_{k=1}^3\int_\Omega |F_k|^2\partial_3 \varphi,
\end{align*}
where the boundary integral vanishes thanks to $F_j\cdot N = 0$. Now we deal with the boundary integral,
\begin{align*}
\int_\Sigma q v\cdot N
		&=-\sigma \int_\Sigma \overline{\nabla} \cdot\frac{\overline{\nabla}\psi}{\sqrt{1+|\overline{\nabla} \psi |^2}} \partial_t\psi 
		+ \int_\Sigma \kappa(1-\overline{\Delta})^2\psi\partial_t \psi  
		+ \int_\Sigma \kappa(1-\overline{\Delta})\partial_t\psi \partial_t\psi  \\
		&=\sigma\int_\Sigma  \frac{\overline{\nabla}\psi}{\sqrt{1+|\overline{\nabla} \psi |^2}} \cdot \overline{\nabla} \partial_t\psi 
		+ \int_\Sigma \kappa(1-\overline{\Delta})\psi(1-\overline{\Delta}) \partial_t \psi 
		+\int_\Sigma \kappa |\japbr\partial_t\psi|^2 \\
		&=\frac{1}{2}\frac{d}{dt}\int_\Sigma \bigg(\sigma\sqrt{1+|\overline{\nabla}\psi|^2}+\kappa|(1-\overline{\Delta})\psi |^2\bigg)
		+\frac{d}{dt}\int_0^t\int_\Sigma \kappa |\japbr\partial_t\psi|^2.
	\end{align*}\\
	Hence, we have the following conserved energy
	\begin{align}
		E_0^\kappa (t)= \frac{1}{2}\int_\Omega \bigg(|v|^2
		+\sum_{k=1}^3|F_k|^2\bigg)\partial_3 \varphi
		+\frac{1}{2}\int_\Sigma \bigg(\sigma\sqrt{1+|\overline{\nabla}\psi|^2}+\kappa|(1-\overline{\Delta})\psi |^2\bigg)
		+\int_0^t\int_\Sigma \kappa |\japbr\partial_t\psi|^2.
	\end{align}

\subsection{Div-curl estimate}
For the higher-order energy estimate, we use the div-curl estimates to convert the normal derivatives to the tangential derivative and curl. 
\begin{lemma}[The Hodge-type elliptic estimate] \label{hodge}
Let $X$ be a smooth vector field and $s\geq 1$. Then 
		\begin{equation}
			||X||_s^2\lesssim C\big(|\psi|_{s},|\overline{\partial}\psi|_{W^{1,\infty}}\big)\big(||\nabla^\varphi\cdot X||_{s-1}^2+||\nabla^\varphi\times X||_{s-1}^2+||\overline{\partial}^s X||_0^2+||X||_0^2\big).\\
		\end{equation}
		For the proof, we refer to \cite{ginsberg2020local}.
\end{lemma}
Applying lemma \ref{hodge} to $||\partial_t^k v||_{4-l}^2$ and $||\partial_t^lF_j||_{4-k}^2$ for $l=0,1,2,3$, we  get 
\begin{align}
		||\partial_t^l v||_{4-l}^2&\lesssim C\big(|\psi|_{4-l},|\overline{\partial}\psi|_{W^{1,\infty}}\big)\big(||\nabla^\varphi\times \partial_t^l v||_{3-l}^2+||\overline{\partial}^{4-l} \partial_t^l v||_0^2+|| \partial_t^lv||_0^2\big),\\
		||\partial_t^lF_j||_{4-l}^2&\lesssim C\big(|\psi|_{4-l},|\overline{\partial}\psi|_{W^{1,\infty}}\big)\big(||\nabla^\varphi\times \partial_t^lF_j||_{3-l}^2+||\overline{\partial}^{4-l} \partial_t^l F_j||_0^2+||\partial_t^lF_j||_0^2\big).
	\end{align}
	We have studied the $L^2$ part in section \ref{L2}. 	For the vorticity part, we take $\nabla^\varphi\times $ to the evolution equations of $v$ and $F_j$ to obtain the evolution equation of the vorticity terms
		\begin{align}
		D_t^{ \varphi}(\partial^{3}\nabla^{\varphi}\times v)&= ({F}_k\cdot\nabla^{ \varphi})(\partial^3\nabla^{\varphi}\times F_k)+ R_1,\\
		D_t^{ \varphi}(\partial^{3}\nabla^{\varphi}\times F_j)&= ({F}_j\cdot\nabla^{ \varphi})(\partial^3\nabla^{\varphi}\times v)+ R_2,
	\end{align}
	where $R_1= \partial^3[\nabla^{\varphi}\times,D_t^{\varphi}]v+\partial^3 [\nabla^{\varphi}\times,( { F}_k\cdot\nabla^{\varphi})]F_k+[\partial^3,D_t^{\varphi}]\nabla^{\varphi}\times v+[\partial^3,({ F}_k\cdot\nabla^{ \varphi})] (\nabla^{\varphi}\times F_k)$ and $R_2 = \partial^3[\nabla^{\varphi}\times,D_t^{\varphi}]F_j+\partial^3 [\nabla^{\varphi}\times,( { F}_j\cdot\nabla^{\varphi})]v+[\partial^3,D_t^{\varphi}]\nabla^{\varphi}\times F_j+[\partial^3,({ F}_j\cdot\nabla^{ \varphi})] (\nabla^{\varphi}\times v)$.
	We test the first equation with  $(\partial^{3}\nabla^{\varphi}\times v)$ and the second equation with $(\partial^{3}\nabla^{\varphi}\times F_j)$ to get
	\begin{align}
		\frac{1}{2}\frac{d}{dt}\int_\Omega |\partial^3\nabla^{\varphi}\times v|^2\partial_3\varphi &= \int_\Omega ({F}_k\cdot\nabla^{ \varphi})(\partial^3\nabla^{\varphi}\times F_k) (\partial^3\nabla^{\varphi}\times v)\partial_3 \varphi+\int_\Omega (\partial^3\nabla^{\varphi}\times v)R_1 \partial_3\varphi\\
		\frac{1}{2}\frac{d}{dt}\int_\Omega |\partial^{3}\nabla^{\varphi}\times F_j|^2\partial_3\varphi &= \int_\Omega ({F}_k\cdot\nabla^{ \varphi})(\partial^3\nabla^{\varphi}\times v) (\partial^3\nabla^{\varphi}\times F_k)\partial_3 \varphi+\int_\Omega (\partial^{3}\nabla^{\varphi}\times F_j)R_1 \partial_3\varphi
	\end{align}
	The top order terms $\int_\Omega ({F}_k\cdot\nabla^{ \varphi})(\partial^3\nabla^{\varphi}\times F_k) (\partial^3\nabla^{\varphi}\times v)\partial_3 \varphi$ and $\int_\Omega ({F}_k\cdot\nabla^{ \varphi})(\partial^3\nabla^{\varphi}\times v) (\partial^3\nabla^{\varphi}\times F_k)\partial_3 \varphi$ cancels each other after integrating $({F}_k\cdot\nabla^{ \varphi})$ by parts and the commutator terms can be directly controlled. It follows that
	\begin{align}
		\frac{1}{2}\frac{d}{dt}\big(||\partial^3\nabla^{\varphi}\times v||_0^2+||\partial^{3}\nabla^{\varphi}\times F_k||_0^2 \big)\lesssim P(||v||_{4},|| F_k||_{4}, |\psi|_5,|\partial_t \psi|_4).
	\end{align}
	We can apply parallel arguments to treat the case with $\partial_t^l\partial^{3-l}$ to get
		\begin{align}
		\frac{1}{2}\frac{d}{dt}\big(||\nabla^{\varphi}\times \partial_t^lv||_{3-l}^2+||\nabla^{\varphi}\times \partial_t^lF_k||_{3-l}^2 \big)\lesssim P(||\partial_t^lv||_{4-l}, ||\partial_t^l F_j||_{4-l}, |\partial_t^l \psi |_{5-l},|\partial_t^4 \psi |_{1} )
	\end{align}

	\subsection{Tangential estimate}\label{prioritangential}
	The higher-order tangential estimate is the key part of the div-curl analysis. In the priori energy estimate \cite{xu2024priori}, we introduced Alinhac's good unknown(AGU) to avoid the appearance of the top order of $\psi$ and get a uniform-in-$\sigma$ energy estimate. In this paper, $\sigma$-uniform estimate is no longer our concern, but we still introduce AGU for neatness and consistency. 
	\subsubsection{Reformulation in Alinhac's good unknown}
	Let $D^\alpha= \partial_t^{\alpha_0}\partial_1^{\alpha_1}\partial_2^{\alpha_2}$ for any multi-index $\alpha$ with $|\alpha|\leq4$, generic function $f$ and $i=1,2,3$. Invoking the definition of $\partial^\varphi$, we have
	\begin{align}
		D^\alpha \partial_i^\varphi f
		=\partial^\varphi_i(D^\alpha f-D^\alpha \varphi\partial_3^\varphi f) +\mathcal{C}_i(f),
	\end{align}
	where, for $\tau =1,2, |\beta|=1$,
	\begin{align}
		\mathcal{C}_\tau (f)&= D^\alpha \varphi\partial_\tau^\varphi\partial_3^\varphi f -[D^\alpha,\frac{\partial_\tau \varphi}{\partial_3 \varphi },\partial_3 f]-\partial_3 f[D^\alpha,\partial_\tau \varphi,\frac{1}{\partial_3 \varphi}]+\partial_3 f \partial _\tau \varphi[D^{\alpha-\beta},\frac{1}{(\partial_3 \varphi)^2}]D^\beta \partial_3 \varphi,\\
		\mathcal{C}_3(f)&=D^\alpha \varphi (\partial_3 ^\varphi)^2 f +[D^\alpha,\frac{1}{\partial_3 \varphi},\partial_3 f]- \partial_3 f [D^{\alpha-\beta},\frac{1}{(\partial_3 \varphi)^2}]D^\beta \partial_3 \varphi,
	\end{align}
	and 
	\begin{align}
		D^\alpha D_t^\varphi f
		=D_t^\varphi (D^\alpha f- D^\alpha \varphi\partial_3 ^\varphi f)+\mathcal{D}(f)\label{D},
	\end{align}
	where 
	\begin{align}
			\mathcal{D}(f)&=D^\alpha \varphi D_t^\varphi \partial_3^\varphi f +[D^\alpha, \overline{v}]\cdot \overline{\nabla}f+[D^\alpha, \frac{1}{\partial_3 \varphi}(v\cdot \mathbf N-\partial_t \varphi),\partial_3 f]+[D^\alpha, \frac{1}{\partial_3 \varphi},v\cdot \mathbf{N}-\partial_t \varphi]\partial_3 f\nonumber\\
			&\quad-(v\cdot \mathbf{N} -\partial_t \varphi)\partial_3 f[D^{\alpha-\beta},\frac{1}{(\partial_3 \varphi)^2}]D^\beta \partial_3 \varphi + \frac{1}{\partial_3 \varphi}\partial_3 f[D^\alpha,v]\mathbf{N}.
	\end{align}
Here $D^\alpha f-D^\alpha \varphi\partial_3^\varphi f$ is the AGU of $f$. It directly follows from the definition that we can obtain the estimate of $f$ once we can control its AGU,
\begin{align}
	||D^\alpha f||_0\leq ||D^\alpha f-D^\alpha \varphi\partial_3^\varphi f||_0 + ||\partial_3^\varphi f||_0||D^\alpha \varphi||_0\label{AGU}.
	\end{align}
 Let $\mathbf{V}, \mathbf{Q}, \mathbf{F}_k$ denotes the AGU of $v,q,F$. 
\begin{align}
	\mathbf{V}=D^\alpha v - D^\alpha \varphi\partial_3^\varphi v,\quad
	\mathbf{F}=D^\alpha F - D^\alpha \varphi\partial_3^\varphi F,\quad
	\mathbf{Q}=D^\alpha q - D^\alpha \varphi\partial_3^\varphi q.
	\end{align}
We write the $D^\alpha$-differentiated system in terms of the AGUs as follows
\begin{align}\label{AGUsys}
\begin{cases}
	D_t^\varphi \mathbf{V}_i + \partial_i^\varphi \mathbf{Q}= (F_k\cdot \nabla^\varphi) \mathbf{F}_{ik}+\mathcal{R}_i^1 \qquad & \text{in } \Omega,\\
	\partial_k^\varphi \mathbf{V}_k=-\mathcal{C}_k(v_k) & \text{in } \Omega,\\
	D_t^\varphi \mathbf{F}_{ij}=(F_j\cdot \nabla^\varphi)\mathbf{V}_i+\mathcal R_{ij}^2  & \text{in } \Omega, \\
	\nabla^\varphi \cdot  F_j=0  & \text{in } \Omega , \\
	\mathbf{Q} = -\sigma D^\alpha\big(\overline{\nabla}\cdot \frac{\overline{\nabla}\psi}{|N|}\big)
	+ \kappa D^\alpha(1-\overline{\Delta})^2\psi + \kappa D^\alpha(1-\overline{\Delta})\partial_t\psi 
	-D^\alpha \psi \partial_3 q   &\text{on }\Sigma , \\
	\partial_t D^\alpha \psi =\mathbf{V}\cdot N-(\overline{v}\cdot \overline{\nabla})  D^\alpha \psi+\mathcal{S}_1  &\text{on }\Sigma ,  \\
	F_j\cdot N =0  &\text{on }\Sigma , \\ 
	\mathbf{V}_3 = \mathbf{F}_{3k}=0  &\text{on }\Sigma_b , 
\end{cases} 
\end{align}
where 
\begin{align}
\mathcal{R}_i^1&=F_{lk}\mathcal{C}_l(F_{ik})+[D^\alpha, F_{lk}] \partial_l^\varphi F_{ik}-\mathcal{D}(v_i)-\mathcal{C}_i(q), \\
\mathcal{R}_{ij}^2&=F_{kj}\mathcal{C}_k(v_i)+[D^\alpha,F_{kj}]\partial^\varphi_k v_i-\mathcal{D}(F_{ij}),\\
\mathcal{S}_1&=D^\alpha \psi \partial _3 v\cdot N -[D^\alpha,\overline{v}\cdot,\overline{\partial}\psi] \label{S1}.
\end{align}
We refer to \cite{xu2024priori} for the detailed derivation.


\subsubsection{Tangential estimates with full spatial derivatives}\label{fullspatial}
We first study the equations with only spatial derivatives, that is, $D^\alpha = \parbar^\alpha$
\begin{prop}
For small $\delta>0$, the following uniform-in-$\kappa$ inequality holds
	\begin{align*}
		&||\mathbf{V}(t)||^2_0+\sum_k||\mathbf{F}_k(t)||^2_0+|\sqrt{\kappa}D^\alpha(1-\overline{\Delta})\psi|_0^2+|\sqrt{\sigma}D^\alpha\overline{\nabla}\psi|^2_0+\int_0^t \int_\Sigma |\sqrt\kappa D^\alpha\japbr\partial_t\psi|_0^2\\
		&\leq \delta E_4^\kappa (t)+\int_0^t P\big(\sigma^{-1},E_4^\kappa(t)\big).
	\end{align*}
\end{prop}

\noindent\textbf{Part 1: Interior structure}\\
We test the first equation of (\ref{AGUsys}) with $\mathbf{V}\partial_3\varphi$ and integrate $\partial_i^\varphi$ by parts to get
\begin{align*}
	\frac{1}{2}\frac{d}{dt}\int_\Omega |\mathbf{V}_i|^2\partial_3 \varphi 
	&= -\int_\Omega\partial_i^\varphi \mathbf{Q} \mathbf{V}_i \partial_3 \varphi + \int_\Omega(F_k\cdot \nabla^\varphi) \mathbf{F}_{ik}\mathbf{V}_i \partial_3 \varphi+\int_\Omega \mathcal{R}_i^1\mathbf{V}_i \partial_3 \varphi \\
	&=-\int_\Omega \mathcal{C}_i(v_i)\mathbf{Q}\partial_3 \varphi-\int_\Sigma \mathbf{Q}\mathbf{V}\cdot N+ \int_\Omega(F_k\cdot \nabla^\varphi) \mathbf{F}_{ik}\mathbf{V}_i \partial_3 \varphi+\int_\Omega \mathcal{R}_i^1\mathbf{V}_i \partial_3 \varphi.
\end{align*}
We invoke the third equation of (\ref{AGUsys}) and integrate $(F_k\cdot \nabla^\varphi)$ by parts to get
\begin{align*}
	\int_\Omega(F_k\cdot \nabla^\varphi) \mathbf{F}_{ik}\mathbf{V}_i \partial_3 \varphi 
	&= -\int_\Omega \mathbf{F}_{ik}D_t^\varphi \mathbf{F}_{ik} \partial_3 \varphi +\int_\Omega \mathbf{F}_{ik} \mathcal{R}^2_{ik} \partial_3 \varphi\\
	&= -\frac{1}{2}\frac{d}{dt} \int_\Omega \sum _{k}|\mathbf{F}_{ik}|^2 \partial_3 \varphi +\int_\Omega \mathbf{F}_{ik} \mathcal{R}^2_{ik} \partial_3 \varphi.
	\end{align*}
	Hence, we have
	\begin{align}
		\frac{1}{2}\frac{d}{dt}\int_\Omega \bigg(|\mathbf{V}_i|^2+\sum_k|\mathbf{F}_{ik}|^2\bigg) \partial_3 \varphi 
		=-\int_\Sigma \mathbf{Q}\mathbf{V}\cdot N 
		-\int_\Omega \mathcal{C}_i(v_i)\mathbf{Q}\partial_3 \varphi
		+\int_\Omega \mathcal{R}_i^1\mathbf{V}_i \partial_3 \varphi
		+\int_\Omega \mathbf{F}_{ik} \mathcal{R}^2_{ik} \partial_3 \varphi.
	\end{align}
	Note that in the case with only spatial derivatives, the commutators can be directly controlled
	\begin{align}
		||\mathcal{C}_i(f)||_0&\leq P(c_0^{-1},|\psi|_4)||f||_4 \label{Cest},\\
		||\mathcal{D}(f)||_0&\leq P(c_0^{-1},||v||_4,|\psi|_4,|\partial_t \psi|_3)(||f||_4+||\partial_t f||_2) \label{Dest},
	\end{align}
	so that we can easily estimate the remainder terms
	\begin{align}
	-\int_\Omega\mathcal{C}_i(v_i)\mathbf{Q}\partial_3 \varphi 
	&\leq||\mathcal{C}_i(v_i)||_0 ||\mathbf{Q}||_0||\partial_3 \varphi||_\infty, \label{civiq}\\	
	\int_\Omega \mathbf{F}_{ik}\mathcal{R}_{ik}^2\partial_3 \varphi 
	&\lesssim||\partial_3 \varphi||_\infty ||\mathbf{F}||_0\big(||F_{kj}||_\infty||\mathcal{C}_k(v_i)||_0 +||F_{kj}||_4||v||_4+||\mathcal{D}(F_{ij})||_0\big) ,\\
	\int_\Omega \mathcal{R}_i^1\mathbf{V}_i \partial_3 \varphi
	&\lesssim||\mathbf{V}_i||_0||\partial_3\varphi||_\infty\big(||F_{ik}||_\infty ||\mathcal{C}_l(F_{ik})||_0+||F_{lk}||_4^2+||\mathcal{D}(v_i)||_0+||\mathcal{C}_i(q)||_0\big).
	\end{align}	
	\\\\
	\textbf{Part 2: Boundary Structure}\\
	The boundary integral is way more complicated than the interior one. Invoking the higher-order kinematic boundary condition, we get 
	\begin{align}\label{bdstruc}
		-\int_\Sigma \mathbf{Q} (\mathbf{V}\cdot N) = -\int_\Sigma \mathbf{Q}\partial_t D^\alpha \psi 
		-\int_\Sigma \mathbf{Q}(\overline{v}\cdot \overline{\nabla})  D^\alpha \psi 
		+\int_\Sigma \mathbf{Q}\mathcal{S}_1.
	\end{align}
	The first term gives 
	\begin{align*}
		-\int_\Sigma \mathbf{Q}\partial_t D^\alpha \psi 
		&= \int_\Sigma \bigg( \sigma  D^\alpha\big(\overline{\nabla}\cdot \frac{\overline{\nabla}\psi}{|N|}\big)
	- \kappa D^\alpha(1-\overline{\Delta})^2\psi - \kappa D^\alpha(1-\overline{\Delta})\partial_t\psi \bigg)\partial_t D^\alpha \psi
	+\int_\Sigma D^\alpha \psi \partial_3 q \partial_t D^\alpha \psi\\
	&=: ST + RT.
	\end{align*}
	The ST term contributes to the boundary regularity,
	\begin{align*}
		ST= \sigma \int_\Sigma  D^\alpha\big(\overline{\nabla}\cdot \frac{\overline{\nabla}\psi}{|N|}\big)D^\alpha \psi
	-  -\frac{1}{2}\frac{d}{dt} \int_\Sigma \kappa|D^\alpha(1-\overline{\Delta})\psi|^2 	
	- \frac{d}{dt} \int_0^t \int_\Sigma \kappa |D^\alpha\japbr\partial_t\psi|^2.
	\end{align*}
	There is a symmetric structure in the first term, 
	\begin{align*}
		\sigma \int_\Sigma  D^\alpha\big(\overline{\nabla}\cdot \frac{\overline{\nabla}\psi}{|N|}\big)D^\alpha \psi 
		& =-\int_\Sigma\sigma D^\alpha\big(\frac{\overline{\nabla}\psi}{|N|}\big)\cdot \partial_t \overline{\nabla}D^\alpha \psi \\
		&= \int_\Sigma \sigma \big(\frac{D^\alpha\overline{\nabla}\psi}{|N|}- \frac{\overline{\nabla}\psi\cdot D^\alpha \overline{\nabla}\psi}{|N|^3}\overline{\nabla}\psi \big)\cdot \partial_t \overline{\nabla}D^\alpha \psi + ST_1^R\\
		&=-\frac{\sigma}{2}\frac{d}{dt}\int_\Sigma \big( \frac{|D^\alpha \overline{\nabla }\psi|^2}{|N|}
		-\frac{|\overline{\nabla}\psi\cdot D^\alpha \overline{\nabla}\psi|^2}{|N|^3} \big) +ST_1^R +ST_2^R,
	\end{align*}
	where the commutators
	\begin{align*}
	ST_1^R&= -\int_\Sigma \sigma\big([D^\alpha,\overline{\nabla}\psi,\frac{1}{|N|}] -\overline{\nabla}\psi [D^{\alpha-\beta},\frac{\overline{\nabla}\psi}{|N|^3}]\cdot D^\beta \overline{\nabla}\psi \big)\cdot \partial_t \overline{\nabla}D^\alpha \psi,\\
		ST_2^R &= \frac{\sigma}{2}\int_\Sigma \partial_t(\frac{1}{|N|})|D^\alpha\overline{\nabla}\psi|^2-\partial_t(\frac{1}{|N|^3})|\overline{\nabla}\psi\cdot D^\alpha \overline{\nabla}\psi|^2.
	\end{align*}
	$ST_2^R$ can be directly controlled and $ST_1^R$ can be controlled by integrating $\overline{\nabla}$ by parts
	\begin{align}
		ST_1^R&\leq P(c_0^{-1},|\psi|_{W^{3,\infty}})(1+|\sqrt{\sigma}\overline{\nabla}D^\alpha \psi|_0|\sqrt{\sigma}\partial_t D^\alpha \psi|_0),\label{str1}\\
		ST_2^R&\leq P(c_0^{-1})|\overline{\nabla}\partial_t \psi |_\infty |\sqrt{\sigma}D^\alpha \overline{\nabla}\psi|_0^2(1+|\overline{\nabla}\psi|_\infty).
	\end{align}
 By Cauchy's inequality, 
\begin{align*}
		\frac{1}{2}\int _\Sigma  \frac{\sigma|D^\alpha \overline{\nabla}\psi|^2 }{|N|^3}
		\leq \frac{\sigma}{2} \int _\Sigma \big(\frac{|D^\alpha \overline{\nabla }\psi|^2}{|N|}-\frac{|\overline{\nabla}\psi\cdot D^\alpha \overline{\nabla}\psi|^2}{|N|^3}\big).
	\end{align*}
	Hence the ST term contributes to the following energy 
	\begin{align}
		\frac{1}{2}\int_\Sigma \kappa|D^\alpha(1-\overline{\Delta})\psi|^2 	+  \frac{1}{2}\int _\Sigma  \frac{\sigma|D^\alpha \overline{\nabla}\psi|^2 }{|N|^3}+\int_0^t \int_\Sigma \kappa |D^\alpha\japbr\partial_t\psi|^2.
	\end{align}
	Since we do not require the estimate to be $\sigma$-uniform, we can directly control the RT term
	\begin{align}
		RT&=\int_\Sigma (\partial_3 q)D^\alpha \psi \partial_t D^\alpha \psi \leq |\partial_3 q|_\infty	|D^\alpha \psi|_0|\partial_t D^\alpha \psi|_0.\label{rt}
	\end{align}
	
	Now we deal with the second term of (\ref{bdstruc}),
	\begin{align} \label{bd2}
		-\int_\Sigma \mathbf{Q}(\overline{v}\cdot \overline{\nabla})  D^\alpha \psi 
		=  \int_\Sigma  \sigma D^\alpha\big(\overline{\nabla}\cdot \frac{\overline{\nabla}\psi}{|N|}\big)(\overline{v}\cdot \overline{\nabla})  D^\alpha \psi	
		- \int_\Sigma \kappa D^\alpha(1-\overline{\Delta})^2\psi(\overline{v}\cdot \overline{\nabla})  D^\alpha \psi  \nonumber\\
	- \int_\Sigma \kappa D^\alpha(1-\overline{\Delta})\partial_t\psi(\overline{v}\cdot \overline{\nabla})  D^\alpha \psi
	+\int_\Sigma D^\alpha \psi \partial_3 q (\overline{v}\cdot \overline{\nabla})  D^\alpha \psi.
	\end{align}
	The second and the third term can be controlled by the $\sqrt\kappa$-weighted energy. Integrating $1-\overline{\Delta}$ by parts, the second term of (\ref{bd2}) gives
	\begin{align*}
		&- \int_\Sigma \kappa D^\alpha(1-\overline{\Delta})^2\psi(\overline{v}\cdot \overline{\nabla})  D^\alpha \psi\\
		&=  - \int_\Sigma \kappa D^\alpha(1-\overline{\Delta})\psi(\overline{v}\cdot \overline{\nabla})  (1-\overline{\Delta})D^\alpha \psi 
		- \int_\Sigma \kappa D^\alpha(1-\overline{\Delta})\psi[(1-\overline{\Delta}),\overline{v}\cdot\overline{\nabla}] D^\alpha\psi,
	\end{align*}
	where the second term can be directly controlled by $|\sqrt\kappa D^\alpha(1-\overline{\Delta})\psi|_0|v|_{W^{2,\infty}}|\sqrt{\kappa}D^\alpha \psi|_2$ and the first term can be controlled by integrating $\overline{v} \cdot\overline{\nabla}$ by parts and symmetry,
	\begin{align*}
		&- 2\int_\Sigma \kappa D^\alpha(1-\overline{\Delta})\psi(\overline{v}\cdot \overline{\nabla})  (1-\overline{\Delta})D^\alpha \psi \\
 		&=\int_\Sigma \kappa D^\alpha(1-\overline{\Delta})\psi (\overline{\nabla}\cdot\overline{v} ) (1-\overline{\Delta})D^\alpha \psi 
 		\leq |\sqrt\kappa D^\alpha(1-\overline{\Delta})\psi|^2_0|v|_{W^{1,\infty}}.
	\end{align*}
	The third term of (\ref{bd2}) under time integral can be handled by integrating $\japbr$ by parts and Cauchy's inequality,
	\begin{align*}
		&- \int_0^t \int_\Sigma \kappa D^\alpha(1-\overline{\Delta})\partial_t\psi(\overline{v}\cdot \overline{\nabla})  D^\alpha \psi  \\
		&= - \int_0^t \int_\Sigma \kappa D^\alpha\japbr\partial_t\psi\japbr\big((\overline{v}\cdot \overline{\nabla})  D^\alpha \psi \big) 
		\leq \delta\int_0^t\int_\Sigma |\sqrt\kappa D^\alpha\japbr\partial_t\psi|^2  
		+\frac{1}{4\delta}\int_0^t|v|^2_{W^{1,\infty}}|\sqrt\kappa D^\alpha\psi|_2^2,
	\end{align*}
	where the first term can be absorbed to $E_4^\kappa(t)$ if we take $\delta$ small enough. We can find a similar symmetric structure in the first term of (\ref{bd2}). Let $|\beta|=1$. We expand $D^\alpha(\frac{\overline{\nabla}\psi}{|N|})$ and integrate $\overline{\nabla}$ by parts to get
	\begin{align}
		\int_\Sigma  \sigma D^\alpha\big(\overline{\nabla}\cdot \frac{\overline{\nabla}\psi}{|N|}\big)(\overline{v}\cdot \overline{\nabla})  D^\alpha \psi
		&=-\sigma \int_\Sigma \big(\frac{D^\alpha \overline{\nabla}\psi}{|N|}		-\frac{\overline{\nabla}\psi\cdot D^\alpha \overline{\nabla}\psi}{|N|^3}\overline{\nabla}\psi\big) \cdot (\overline{v}\cdot \overline{\nabla})D^\alpha \overline{\nabla}\psi \nonumber \\
	&\quad-\sigma \int_\Sigma \bigg(\big(\frac{D^\alpha \overline{\nabla}\psi}{|N|}		-\frac{\overline{\nabla}\psi\cdot D^\alpha \overline{\nabla}\psi}{|N|^3}\overline{\nabla}\psi\big) \cdot \overline{\nabla} \bigg)\overline{v}\cdot \overline{\nabla}D^\alpha \psi\nonumber\\
	&\quad +\sigma\int_\Sigma \overline{\nabla}\cdot \big([D^\alpha,\overline{\nabla}\psi,\frac{1}{|N|}] -\overline{\nabla}\psi [D^{\alpha-\beta},\frac{\overline{\nabla}\psi}{|N|^3}]\cdot D^\beta \overline{\nabla}\psi\big)(\overline{v}\cdot \overline{\nabla})D^\alpha \psi, \nonumber
		\end{align}
	where the second and the third term can be directly controlled, and the first term can be controlled by integrating $(\overline{v}\cdot \overline{\nabla})$ by parts to get a symmetric structure. Hence,
	\begin{align*}
		\int_\Sigma  \sigma D^\alpha\big(\overline{\nabla}\cdot \frac{\overline{\nabla}\psi}{|N|}\big)(\overline{v}\cdot \overline{\nabla})  D^\alpha \psi \leq P(|\overline{\nabla}\psi|_{W^{1,\infty}})|v|_{W^{1,\infty}} |\sqrt{\sigma}D^\alpha \overline{\nabla}\psi|_0^2.
	\end{align*} 
	The last term of (\ref{bd2}) can be directly controlled 
	\begin{align*}
	\int_\Sigma D^\alpha \psi \partial_3 q (\overline{v}\cdot \overline{\nabla})  D^\alpha \psi \leq \sigma^{-1}|\partial_3q|_\infty|v|_\infty|\sqrt\sigma D^\alpha\psi|_0|\sqrt\sigma D^\alpha\overline{\nabla}\psi|_0	.
	\end{align*}

For the last term of (\ref{bdstruc}), $\mathcal{S}_1$ consists of lower order terms except $\parbar^3 v $, but it can still take $\frac{1}{2}$ of spatial derivative to reach the top order,
	\begin{align*}
		\int_\Sigma \mathbf{Q}\mathcal{S}_1 
	= \int_\Sigma D^\alpha q \mathcal{S}_1 - \int_\Sigma D^\alpha \psi \partial_3 q\mathcal{S}_1.
	\end{align*}
	The second term can be directly controlled,
	 \begin{align*}
	 	\int_\Sigma D^\alpha \psi \partial_3 q\mathcal{S}_1
	 	&\leq |D^\alpha\psi|_0|\partial_3 q|_\infty|\mathcal S_1 |_0\\
	 	&\leq|D^\alpha\psi|_0|\partial_3 q|_\infty\big(P(|\parbar\psi|_\infty)|D^\alpha\psi|_0|v|_{W^{1,\infty}} +|v|_{W^{1,\infty}}|\parbar\psi|_3+2|\parbar\psi|_{W^{2,\infty}}|v|_3\big).
	 \end{align*}
	 For the first term, since $D^\alpha$ contains spatial derivative, we integrate $\parbar^\frac{1}{2}$ by parts to get
	 \begin{align}
	 	\int_\Sigma D^\alpha q \mathcal{S}_1 \leq |D^\alpha q|_{-\frac{1}{2}}|\mathcal S_1|_\frac{1}{2} \leq |q|_{3.5}\big(P(|\parbar \psi|_{W^{1,\infty}})|D^\alpha\psi|_1|v|_{W^{2,\infty}}  +|v|_{W^{2,\infty}}|\parbar\psi|_4+2|\parbar\psi|_{W^{3,\infty}}|v|_{3.5}\big)\label{fulltq}.
	 \end{align}
	 
\subsubsection{Tangential estimates with full time derivatives}\label{fulltime}
The arguments in section (\ref{fullspatial}) still apply when there is at least one spatial derivative in $D^\alpha$. However, there are extra difficulties with full-time derivatives, that is, $D^\alpha = \partial_t^4$. \\\\
\textbf{Difficulties}
 \begin{enumerate}[label=(\roman*)]
 	\item  In the estimates (\ref{str1}), (\ref{rt}) and $||\mathcal C_i(v_i)||_0$, there are appearance of $|\sqrt\sigma \partial_t D^\alpha\psi|_0$ or $|D^\alpha\psi|_0$, which cannot be directly controlled in the case $D^\alpha =\partial_t^4$, since our $\sigma$-weighted energy requires $\psi$ to be attached with $\overline{\nabla}$.
 	\item For the remainder $\int_\Omega\mathcal{C}_i(v_i)\mathbf Q \partial_3 \varphi$, we no longer have direct control for $\mathbf{Q} = \partial_t^4 q -\partial_t^4\psi \partial_3 q$ due to $\partial_t^4 q$.
 	\item In (\ref{fulltq}), we cannot integrate $\partial^\frac{1}{2}$ by parts when there are only time derivatives in $D^\alpha$.
 \end{enumerate}
For (i), to deal with the trouble terms $\partial_t\psi,\partial_t^5\psi$, we invoke the kinematic boundary condition to trade 1 time derivative for 1 spatial derivative
\begin{align}
\partial_t^4\psi
&=\partial_t ^3 (v\cdot N)=-(\overline{v}\cdot \overline{\partial})\partial_t ^3 \psi+\partial_t^3v \cdot N -[\partial_t^3,\overline{v},\overline{\partial}\psi ]\label{t4psi},\\
\partial_t^5\psi
& = \partial_t^4(v\cdot N)=
-(\overline{v}\cdot \overline{\partial})\partial_t ^4 \psi+\partial_t^4 v \cdot N -[\partial_t^4,\overline{v},\overline{\partial}\psi ].\label{t5psi}
\end{align}
we obtained a direct control for $\partial_t^4\psi$,
\begin{align}
|\partial_t^4\psi|_0\leq |v|_\infty|\partial_t^3 \overline{\nabla}\psi|_0+|\overline{\nabla}\psi|_\infty |\partial_t^3 v|_0+|\partial_t v|_\infty|\partial_t^2\overline{\nabla}\psi|_0+|\partial_t^2 v|_0|\partial_t\overline{\nabla}\psi|_\infty.
\end{align}
but for $\partial_t^5\psi$, the trade produced $\partial_t^4v$, which we cannot control on the boundary and needs further treatment. For (ii) and (iii), we seek a cancellation structure for the trouble terms.

\begin{prop}
For small $\delta>0$, the following uniform-in-$\kappa$ inequality holds,
	\begin{align*}
		&||\mathbf{V}(t)||^2_0+\sum_k||\mathbf{F}_k(t)||^2_0+|\sqrt{\kappa}\partial_t^4(1-\overline{\Delta})\psi|_0^2+|\sqrt{\sigma}\partial_t^4\overline{\nabla}\psi|^2_0+\int_0^t \int_\Sigma |\sqrt\kappa \partial_t^4 \japbr\partial_t\psi|_0^2\\
		&\leq \delta E_4^\kappa (t)+\int_0^t P\big(\sigma^{-1},E_4^\kappa(t)\big).
	\end{align*}
\end{prop}

In this section, we skip the estimates of the terms that can be handled similarly as in the section \ref{fullspatial} and focus only on the difficult terms that arise from full-time derivatives mentioned at the beginning of the section.\\\\
\textbf{Part 1: Interior structure}\\
Let 
\begin{align}
	\mathbf{V}=\partial_t^4 v - \partial_t^4 \varphi\partial_3^\varphi v,\quad
	\mathbf{F}=\partial_t^4 - \partial_t^4 \varphi\partial_3^\varphi F,\quad
	\mathbf{Q}=\partial_t^4 - \partial_t^4 \varphi\partial_3^\varphi q.
	\end{align}
We have the same interior structure
\begin{align}\label{fulltimeinterior}
		\frac{1}{2}\frac{d}{dt}\int_\Omega \bigg(|\mathbf{V}_i|^2+\sum_k|\mathbf{F}_{ik}|^2\bigg) \partial_3 \varphi 
		=-\int_\Sigma \mathbf{Q}\mathbf{V}\cdot N 
		-\int_\Omega \mathcal{C}_i(v_i)\mathbf{Q}\partial_3 \varphi
		+\int_\Omega \mathcal{R}_i^1\mathbf{V}_i \partial_3 \varphi
		+\int_\Omega \mathbf{F}_{ik} \mathcal{R}^2_{ik} \partial_3 \varphi.
	\end{align} 
	where 
\begin{align}
\mathcal{R}_i^1&=F_{lk}\mathcal{C}_l(F_{ik})+[\partial_t^4, F_{lk}] \partial_l^\varphi F_{ik}-\mathcal{D}(v_i)-\mathcal{C}_i(q), \\
\mathcal{R}_{ij}^2&=F_{kj}\mathcal{C}_k(v_i)+[\partial_t^4,F_{kj}]\partial^\varphi_k v_i-\mathcal{D}(F_{ij}),
\end{align}
and the commutators, for $\tau =1,2$,
	\begin{align}
		\mathcal{C}_\tau (f)&= \partial_t^4 \varphi\partial_\tau^\varphi\partial_3^\varphi f -[\partial_t^4 ,\frac{\partial_\tau \varphi}{\partial_3 \varphi },\partial_3 f]-\partial_3 f[\partial_t^4 ,\partial_\tau \varphi,\frac{1}{\partial_3 \varphi}]+\partial_3 f \partial _\tau \varphi[\partial_t^3,\frac{1}{(\partial_3 \varphi)^2}]\partial_t \partial_3 \varphi,\\
		\mathcal{C}_3(f)&=\partial_t^4  \varphi (\partial_3 ^\varphi)^2 f +[\partial_t^4 ,\frac{1}{\partial_3 \varphi},\partial_3 f]- \partial_3 f [\partial_t^3,\frac{1}{(\partial_3 \varphi)^2}]\partial_t \partial_3 \varphi,\\
			\mathcal{D}(f)&=\partial_t^4  \varphi D_t^\varphi \partial_3^\varphi f +[\partial_t^4 , \overline{v}]\cdot \overline{\nabla}f+[\partial_t^4 , \frac{1}{\partial_3 \varphi}(v\cdot \mathbf N-\partial_t \varphi),\partial_3 f]+[\partial_t^4 , \frac{1}{\partial_3 \varphi},v\cdot \mathbf{N}-\partial_t \varphi]\partial_3 f\nonumber\\
			&\quad-(v\cdot \mathbf{N} -\partial_t \varphi)\partial_3 f[\partial_t^3,\frac{1}{(\partial_3 \varphi)^2}]\partial_t \partial_3 \varphi + \frac{1}{\partial_3 \varphi}\partial_3 f[\partial_t^4 ,v]\mathbf{N}.
	\end{align}
	As expected, we can still direct control for the commutators and remainders,
	
\begin{align*}
||\mathcal{C}_i(g)||_0&\leq P\bigg(c_0^{-1},|\overline{\nabla}\psi|_\infty, \sum_{k=1}^3 |\overline{\nabla}\partial_t^k\psi|_{3-k}, |\partial_t^4\psi |_0\bigg) \cdot \bigg(||\partial g||_\infty +||\partial^2 g||_\infty + \sum _{k=1}^3 ||\partial_t^k g||_{4-k}\bigg),\\
||\mathcal{D}(g)||_0&\leq P\bigg(c_0^{-1},\sum_{k=0}^3||\partial_t^k v||_0,|\overline{\nabla}\psi|_\infty, \sum_{k=1}^3 |\overline{\nabla}\partial_t^k\psi|_{3-k}, |\partial_t^4\psi |_0\bigg) \cdot \bigg(||\partial g||_\infty +||\partial^2 g||_\infty + \sum _{k=1}^3 ||\partial_t^k g||_{4-k}\bigg),\\
||\mathcal{R}_i^1||_0 &\leq P\bigg(\sum_{k=0}^4||\partial_t ^k F||_{4-k}\bigg)\cdot \bigg(  ||\mathcal{C}_l(F_{ik})||_0 + ||\mathcal{D}(v_i)||_0 +\mathcal{C}_i(q)||_0\bigg),\\
||\mathcal{R}_{ij}^2||_0 &\leq P\bigg(\sum_{k=0}^4||\partial_t ^k F||_{4-k}, \sum_{k=0}^3||\partial_t^k v||_{4-k}\bigg)\cdot \bigg(  ||\mathcal{C}_l(v_i)||_0 + ||\mathcal{D}(F_{ij})||_0\bigg).
\end{align*}
Hence, the third and last terms of (\ref{fulltimeinterior}) can be directly controlled
\begin{align*}
\int_\Omega \mathbf{F}_{ik} \mathcal{R}^2_{ik} \partial_3 \varphi&\leq ||\partial_3 \varphi||_\infty||\mathbf{F}_{ik}||_0||\mathcal{R}_{ik}^2||_0,\\
\int_\Omega \mathcal{R}_i^1\mathbf{V}_i \partial_3 \varphi 	&\leq ||\partial_3 \varphi ||_\infty ||\mathbf{V}_i||_0||\mathcal{R}_i^1||_0.
\end{align*}
As we mentioned, the only interior trouble term is
\begin{align}
-\int _\Omega \mathcal{C}_i(v_i)\mathbf{Q}\partial_3 \varphi	= -\int_\Omega \partial_t^4q \mathcal{C}_i(v_i) \partial_3\varphi
+\int_\Omega \partial_t^4 \psi \partial_3 q\mathcal{C}_i(v_i)\partial_3 \varphi,
\end{align}
the second term can be directly controlled by $|\partial_t^4 \psi|_0|\partial_3q|_\infty||\mathcal{C}_i(v_i)||_0|\partial_3\varphi|_\infty$, for the first term, we leave it here for the moment and seek a cancellation structure on the boundary. 
\\\\
\textbf{Part 2: Boundary Structure}\\
Invoking the time-differentiated kinematic boundary condition,  
\begin{align}
\partial_t^5\psi = \mathbf{V}\cdot N- \overline{v}\cdot \partial_t^4\overline{\nabla}\psi +\mathcal{S}^*_1,
\end{align}
where 
\begin{align}
\mathcal{S}^*_1=\partial_t^4 \psi \partial_3 v\cdot N-[\partial_t^4, \overline{v},\overline{\partial}\psi].
\end{align}
we have
\begin{align}\label{bdstru2}
	-\int_\Sigma \mathbf{Q}\big(\mathbf{V}\cdot N\big)&= -\int_\Sigma \mathbf{Q} \partial_t^5\psi -\int_\Sigma \mathbf{Q}\big(\overline{v}\cdot \partial_t^4\overline{\nabla}\psi\big) +\int_\Sigma \mathbf{Q}\mathcal{S}^*_1.
\end{align}
As pointed out at the beginning of the section, extra difficulties come from the first and the third term. The second term can be treated by the arguments in section \ref{fullspatial}. By the $\partial^4_t$-differentiated boundary condition
\begin{align}
		\mathbf{Q} = -\sigma \partial^4_t\big(\overline{\nabla}\cdot \frac{\overline{\nabla}\psi}{|N|}\big)
	+ \kappa \partial^4_t(1-\overline{\Delta})^2\psi + \kappa \partial^4_t(1-\overline{\Delta})\partial_t\psi 
	-\partial^4_t \psi \partial_3 q  ,
\end{align}
The first term is split into ST and RT terms,
\begin{align*}
		-\int_\Sigma \mathbf{Q}\partial_t^5 \psi 
		&= \int_\Sigma \bigg( \sigma  \partial_t^4\big(\overline{\nabla}\cdot \frac{\overline{\nabla}\psi}{|N|}\big)
	- \kappa \partial_t^4(1-\overline{\Delta})^2\psi - \kappa\partial_t^4(1-\overline{\Delta})\partial_t\psi \bigg)\partial_t^5 \psi
	+\int_\Sigma \partial_t^4 \psi \partial_3 q \partial_t^5\psi\\
	&=: ST^* + RT^*.
\end{align*}
We expected $ST^*$ to contribute to the boundary regularities by parallel arguments
\begin{align*}
		ST=& -\frac{\sigma}{2}\frac{d}{dt}\int_\Sigma \big( \frac{|\partial_t^4\overline{\nabla }\psi|^2}{|N|}
		-\frac{|\overline{\nabla}\psi\cdot \partial_t^4 \overline{\nabla}\psi|^2}{|N|^3} \big)
	  -\frac{1}{2}\frac{d}{dt} \int_\Sigma \kappa|\partial_t^4(1-\overline{\Delta})\psi|^2 \\
	&- \frac{d}{dt} \int_0^t \int_\Sigma \kappa |\partial_t^4\japbr\partial_t\psi|^2 +ST_1^R +ST_2^R.
	\end{align*}
	where
	\begin{align*}
	ST_1^{R*}&= -\int_\Sigma \sigma\big([\partial_t^4,\overline{\nabla}\psi,\frac{1}{|N|}] -\overline{\nabla}\psi [\partial_t^3,\frac{\overline{\nabla}\psi}{|N|^3}]\cdot \partial_t \overline{\nabla}\psi \big)\cdot \partial_t^5 \overline{\nabla} \psi,\\
		ST_2^{R*} &= \frac{\sigma}{2}\int_\Sigma \partial_t(\frac{1}{|N|})|\partial_t^4\overline{\nabla}\psi|^2-\partial_t(\frac{1}{|N|^3})|\overline{\nabla}\psi\cdot \partial_t^4\overline{\nabla}\psi|^2.
	\end{align*}
The second term is controlled by $P(c_0^{-1})|\overline{\nabla}\partial_t \psi |_\infty |\sqrt{\sigma}\partial_t^4 \overline{\nabla}\psi|_0^2(1+|\overline{\nabla}\psi|_\infty)$. In section \ref{fullspatial}, we integrated $\overline{\nabla}$ by parts to control $ST_1^{R}$ and produced the term $\sqrt{\sigma}\partial_tD^\alpha\psi$, which we cannot control in the case of $D^\alpha= \partial_t^4.$ The straightforward way to solve the problem is that we integrate $\partial_t$ by parts under time integral instead of $\overline{\nabla}.$
\begin{align*}
	&-\int_0^t \int_\Sigma \sigma\big([\partial_t^4,\overline{\nabla}\psi,\frac{1}{|N|}] -\overline{\nabla}\psi [\partial_t^3,\frac{\overline{\nabla}\psi}{|N|^3}]\cdot \partial_t \overline{\nabla}\psi \big)\cdot \partial_t^5 \overline{\nabla} \psi \\
	&\leq P(E(0))+ \int_0^t \int_\Sigma \sigma\partial_t\big([\partial_t^4,\overline{\nabla}\psi,\frac{1}{|N|}] -\overline{\nabla}\psi [\partial_t^3,\frac{\overline{\nabla}\psi}{|N|^3}]\cdot \partial_t \overline{\nabla}\psi \big)\cdot \partial_t^4 \overline{\nabla} \psi \\
	&\quad-  \int_\Sigma \sigma\big([\partial_t^4,\overline{\nabla}\psi,\frac{1}{|N|}] -\overline{\nabla}\psi [\partial_t^3,\frac{\overline{\nabla}\psi}{|N|^3}]\cdot \partial_t \overline{\nabla}\psi \big)\cdot \partial_t^4 \overline{\nabla} \psi,
\end{align*}
where the second term can by bounded by $\int_0^t P(c_0^{-1}, |\partial_t\overline{\nabla}\psi|_\infty, |\partial_t^2\overline{\nabla}\psi|_\infty)(|\sqrt\sigma \partial_t^3\overline{\nabla}\psi|_0+|\sqrt\sigma \partial_t^4\overline{\nabla}\psi|_0)$. For the third term, invoking the Cauchy's inequality
\begin{align*}
  &\int_\Sigma \sigma\big([\partial_t^4,\overline{\nabla}\psi,\frac{1}{|N|}] -\overline{\nabla}\psi [\partial_t^3,\frac{\overline{\nabla}\psi}{|N|^3}]\cdot \partial_t \overline{\nabla}\psi \big)\cdot \partial_t^4 \overline{\nabla} \psi
  \\
  &\leq \delta\int_\Sigma |\partial_t^4\overline{\nabla}\psi|^2 + \frac{1}{4\delta}P(c_0^{-1},|\partial_t \overline{\nabla} \psi|_\infty,|\partial_t^2 \overline{\nabla} \psi|_\infty)|\sqrt\sigma \partial_t^3\overline{\nabla}\psi|_0.
  \end{align*}
Another term that contains $\partial_t^5 $ is $RT^*$, here we invoke the kinematic boundary condition (\ref{t5psi}) to trade a time derivative with a spatial derivative,
\begin{align*}
	RT^* = \int_\Sigma  -\partial_t^4 \psi \partial_3 q(\overline{v}\cdot \overline{\partial})\partial_t ^4 \psi
	+\int_\Sigma \partial_t^4 \psi \partial_3 q\partial_t^4 (v \cdot N) 
	-\int_\Sigma \partial_t^4 \psi \partial_3 q[\partial_t^4,\overline{v},\overline{\partial}\psi ],
\end{align*}
the first term can be controlled by $\sigma^{-\frac{1}{2}}|\partial_t^4\psi|_0|\partial_3 q|_\infty|v|_\infty|\sqrt{\sigma}\partial_t^4\overline{\nabla}\psi|_0$ and the third term can be controlled by $|\partial_t^4\psi|_0|\partial_3 q|_\infty(|\partial_t v|_\infty|\partial_t^3\overline{\nabla}\psi|_0+ |\partial_t^2 v|_\infty|\partial_t^2\overline{\nabla}\psi|_0+|\partial_t^3v|_0|\partial_t\overline{\nabla}\psi|_\infty)$. The second term contains $\partial_t^4 v$, which cannot be controlled on the boundary, we invoke (\ref{t4psi}) and Green's theorem to treat it in the interior 
\begin{align*}
	&\int_\Sigma \big(-(\overline{v}\cdot \overline{\partial})\partial_t ^3 \psi+\partial_t^3v \cdot N -[\partial_t^3,\overline{v},\overline{\partial}\psi ]\big) \partial_3 q\partial_t^4 v \cdot N \\
	& =  \int_\Omega \partial_3\bigg(\big(-(\overline{v}\cdot \overline{\partial})\partial_t ^3 \varphi+\partial_t^3v \cdot \mathbf N -[\partial_t^3,\overline{v},\overline{\partial}\varphi ]\big) \partial_3 q\partial_t^4 v \cdot \mathbf N \bigg) \\
	& = \int_\Omega \partial_3 \partial_t^4 v \cdot \big(-(\overline{v}\cdot \overline{\partial})\partial_t ^3 \varphi+\partial_t^3v \cdot \mathbf N -[\partial_t^3,\overline{v},\overline{\partial}\varphi ]\big) \partial_3 q \mathbf N  \\
	&\quad+\int_\Omega\partial_t^4 v \cdot  \partial_3 \bigg(\big(-(\overline{v}\cdot \overline{\partial})\partial_t ^3 \varphi+\partial_t^3v \cdot \mathbf N -[\partial_t^3,\overline{v},\overline{\partial}\varphi ]\big) \partial_3 q\mathbf N \bigg),
\end{align*}
the second term can be controlled by $P(||v||_{W^{1,\infty}},\sum_{l=1}^3||\partial_t^lv||_1,\sum_{l=0}^2|\partial_t^l\overline{\nabla}\psi|_\infty|\partial_t^3 \overline{\nabla}\psi|_0)||q||_{W^{2,\infty}}|\partial_t^4 v|_0$ and the first term can be controlled by integrating $\partial_t$ by parts under time integral
\begin{align*}
	&\int_0^t\int_\Omega \partial_3 \partial_t^4 v \cdot \big(-(\overline{v}\cdot \overline{\partial})\partial_t ^3 \varphi+\partial_t^3v \cdot \mathbf N -[\partial_t^3,\overline{v},\overline{\partial}\varphi ]\big) \partial_3 q \mathbf N  \\
	&= P(E_4^\kappa(0))-\int_0^t\int_\Omega \partial_3 \partial_t^3 v \cdot \partial_t\bigg(\big(-(\overline{v}\cdot \overline{\partial})\partial_t ^3 \varphi+\partial_t^3v \cdot \mathbf N -[\partial_t^3,\overline{v},\overline{\partial}\varphi ]\big) \partial_3 q \mathbf N \bigg) \\
	&\quad+\int_\Omega \partial_3 \partial_t^3 v \cdot \bigg(\big(-(\overline{v}\cdot \overline{\partial})\partial_t ^3 \varphi+\partial_t^3v \cdot \mathbf N -[\partial_t^3,\overline{v},\overline{\partial}\varphi ]\big) \partial_3 q \mathbf N \bigg),
\end{align*}
where the first term is controlled by $\int_0^t P\bigg(\sum_{k=0}^4||\partial_t^kv||_{4-k},\sum_{k=0}^4|\partial_t^k\overline{\nabla}\psi|_{4-k},||\partial_3 q||_\infty, ||\partial_t\partial_3 q||_\infty \bigg)$ and the top order term in the second term can be singled out by Cauchy's inequality and absorbed into the energy terms,
\begin{align*}
	&\int_\Omega \partial_3 \partial_t^3 v \cdot \bigg(\big(-(\overline{v}\cdot \overline{\partial})\partial_t ^3 \varphi+\partial_t^3v \cdot \mathbf N -[\partial_t^3,\overline{v},\overline{\partial}\varphi ]\big) \partial_3 q \mathbf N \bigg)\\
	&\leq  \delta  ||\partial_3\partial_t^3 v||_0^2 +\frac{1}{4\delta} P\bigg(\sum_{k=0}^3||\partial_t^kv||_{3-k},\sum_{k=0}^3|\partial_t^k\overline{\nabla}\psi|_{3-k}\bigg)||\partial_3q||_\infty	|\overline{\nabla}\psi|_\infty.\label{J2}
\end{align*}
Hence, we have the following estimate
\begin{align}
	 &\int_0^t (ST^*+ RT^*)+\frac{1}{2}\int_\Sigma \frac{|\sqrt{\sigma}\partial_t^4 \overline{\nabla }\psi|^2}{|N|^3}
	  +\frac{1}{2}\int_\Sigma \kappa|\partial_t^4(1-\overline{\Delta})\psi|^2 
	  +  \int_0^t \int_\Sigma \kappa |\partial_t^4\japbr\partial_t\psi|^2\nonumber\\
	  &\leq \delta(||\partial_3\partial_t^3 v||_0^2+ |\partial_t^4\overline{\nabla}\psi|^2_0) +\int_0^tP(\sigma^{-1},E_4^\kappa (t)).
\end{align}
For the last term of (\ref{bdstru2})
\begin{align*}
	\int_\Sigma \mathbf{Q}\mathcal{S}^*_1 = \int_\Sigma\partial_t^4q \mathcal{S}_1^*-\int_\Sigma \partial_t^4\psi\partial_3 q\mathcal{S}_1^*,
\end{align*}
the second term can be controlled by $|\partial_t^4 \psi|_0|\partial_3q|_\infty |\mathcal{S}_1^*|_0$ and $\mathcal{S}_1^*$ consists of only lower order terms. Here we see $\partial_t^4 q$ in the first term that we do not have a direct control,
\begin{align}\label{final}
	\int_\Sigma\partial_t^4q \mathcal{S}_1^* = \int_\Sigma \partial_t^4 q\partial_t^4 \psi \partial _3 v\cdot N 
	-\int_\Sigma \partial_t^4 q[\partial_t^4, \overline{v}\cdot,\overline{\partial}\psi],
\end{align}
In the first term, we see that the terms attached to $\partial _t^4 q$ are of lower order, we invoke the boundary condition for pressure, 
\begin{align*}
	&-\int_\Sigma \partial_t^4 q\partial_t^4 \psi \partial_3 v\cdot N\\
	&=  \int_\Sigma \sigma \partial^4_t\big(\overline{\nabla}\cdot \frac{\overline{\nabla}\psi}{|N|}\big)\partial_t^4 \psi \partial_3 v\cdot N
	- \int_\Sigma\kappa \partial^4_t(1-\overline{\Delta})^2\psi\partial_t^4 \psi \partial_3 v\cdot N
	 - \int_\Sigma\kappa \partial^4_t(1-\overline{\Delta})\partial_t\psi \partial_t^4 \psi \partial_3 v\cdot N,
\end{align*}
these three terms can be controlled by integrating the spatial derivatives by parts.
\begin{align*}
	\int_\Sigma \sigma \partial^4_t\big(\overline{\nabla}\cdot \frac{\overline{\nabla}\psi}{|N|}\big)\partial_t^4 \psi \partial_3 v\cdot N 
	&= \int _\Sigma \sigma \partial_t^4 \big( \frac{\overline{\nabla}\psi}{|N|}\big)\cdot \overline{\nabla}\big(\partial_t ^4 \psi \partial_3 v\cdot N\big)\\
	&\leq 
	P\big(c_0^{-1},|\sqrt{\sigma}\partial_t^4 \overline{\nabla}\psi|_0, |\partial_t ^4 \psi|_0,|\partial_3 \overline{\nabla}v|_\infty,|\partial_3 v|_\infty,|\overline{\nabla}^2\psi|_\infty,|\overline{\nabla}\psi|_\infty \big),\\
	- \int_\Sigma\kappa \partial^4_t(1-\overline{\Delta})^2\psi\partial_t^4 \psi \partial_3 v\cdot N
	&= - \int_\Sigma\kappa \partial^4_t(1-\overline{\Delta})\psi(1-\overline\Delta)(\partial_t^4 \psi \partial_3 v\cdot N)\\
	&\leq P(|\overline{\nabla}\psi|_{W^{2,\infty}},|\sqrt{\kappa}\partial_t^4\psi|_\infty, |\sqrt{\kappa}\partial_t^4\psi|_1,|\sqrt{\kappa}\partial_t^4(1-\overline{\Delta})\psi|_0|,|\partial_3 v|_{W^{1,\infty}}, |v|_3),\\
	- \int_0^t \int_\Sigma\kappa \partial^4_t(1-\overline{\Delta})\partial_t\psi \partial_t^4 \psi \partial_3 v\cdot N 
	&= - \int_0^t \int_\Sigma\kappa \japbr\partial^5_t\psi \japbr( \partial_t^4 \psi \partial_3 v\cdot N)\\
	&\leq \delta \int_0^t \int_\Sigma |\sqrt\kappa \japbr\partial_t^5\psi|^2+\frac{1}{4\delta} \int_0^t \int_\Sigma\kappa |\japbr( \partial_t^4 \psi \partial_3 v\cdot N)|^2\\
	&\leq  \delta \int_0^t \int_\Sigma |\sqrt\kappa \japbr\partial_t^5\psi|^2+\frac{1}{4\delta}\int_0^t |\sqrt\kappa\partial_t^4\psi|_1^2|\partial_3v|_{W^{1,\infty}}^2|\overline{\nabla}\psi|_{W^{1,\infty}}^2.
\end{align*}
For the second term of (\ref{final}), we can single out the term $4 \partial_t^3 \overline v\cdot  \partial_t \overline{\nabla}\psi$, since all the other terms in $[\partial_t^4, \overline{v}\cdot,\overline{\partial}\psi]$ can take one extra derivative to reach the top order, 
\begin{align}\label{final2}
-\int_\Sigma \partial_t^4 q[\partial_t^4, \overline{v}\cdot,\overline{\partial}\psi] 
=-\int_\Sigma 4\partial_t^4 q \partial_t^3 \overline v\cdot  \partial_t \overline{\nabla}\psi 
-\int _\Sigma \partial_t^4q \sum _{k=1}^2 \begin{pmatrix} 4 \\ k \end{pmatrix} \partial_t^kv\cdot \partial_t^{4-k}\overline{\partial}\psi,
\end{align}
we can control the second term by integrating $\partial_t$ by parts under time integral
\begin{align*}
&-\int_0^t\int _\Sigma \partial_t^4q \sum _{k=1}^2 \begin{pmatrix} 4 \\ k \end{pmatrix} \partial_t^kv\cdot \partial_t^{4-k}\overline{\partial}\psi \\
&\leq P(E_4^\kappa(0))+ 	\int_0^t\int _\Sigma \partial_t^3q \sum _{k=1}^2 \begin{pmatrix} 4 \\ k \end{pmatrix} \partial_t\big(\partial_t^kv\cdot \partial_t^{4-k}\overline{\partial}\psi\big)
-\int _\Sigma \partial_t^3q \sum _{k=1}^2 \begin{pmatrix} 4 \\ k \end{pmatrix} \partial_t^kv\cdot \partial_t^{4-k}\overline{\partial}\psi,
\end{align*}
where the second term can be controlled by $\int_0^t|\partial_t^3 q|_0P(\sum_{l=1}^2|\partial_t^lv|_\infty,|\partial_t^3v|_0,|\partial_t^2\overline{\nabla}\psi|_\infty,\sum_{l=3}^4|\partial_t^l\overline{\nabla}\psi|_0)$ and the third term consists of lower terms except $\partial_t^3$ which can be handled by Cauchy's inequality
\begin{align*}
-\int _\Sigma \partial_t^3q \sum _{k=1}^2 \begin{pmatrix} 4 \\ k \end{pmatrix} \partial_t^kv\cdot \partial_t^{4-k}\overline{\partial}\psi \leq \delta |\partial_t^3 q|_0+\frac{1}{4\delta}P(\sum_{l=1}^2|\partial_t^lv|_\infty,\sum_{l=2}^3|\partial_t^l\overline{\nabla}\psi|_0).
\end{align*}
For the first term of (\ref{final2}), the problem is not only that $\partial_t^4 q$ exceeds the top order due to the loss of $\frac{1}{2}$ derivative on the boundary, we do not have the control for $\partial_t^4 q$ as well in the interior. Thus, it does not help to apply Green's theorem to drag it back to the interior as we did for $RT^*$, we can only seek cancellation structures from $-\int_\Omega \partial_t^4 q \mathcal{C}_i(v_i)\partial_3 \varphi$ that we left in the interior estimates. Note that only the following terms need further analysis, 
\begin{align*}
	J_1^*&=-\int_\Omega \partial_t^4q \partial_t^4 \varphi \partial_i^\varphi \partial_3^\varphi v_i \partial_3\varphi, \qquad i=1,2,3,\\
J_2^*&=4\int_\Omega \partial_t^4q \partial_t\big(\frac{\partial_\tau\varphi}{\partial_3 \varphi})\partial_t^3\partial_3 v_\tau \partial_3\varphi, \qquad \tau=1,2,\\
J_3^*&=-4\int_\Omega \partial_t^4q \partial_t\big(\frac{1}{\partial_3\varphi}\big)\partial_t^3 \partial_3 v_3 \partial_3\varphi,	
\end{align*}
the rest are contributed by the lower terms of $\mathcal{C}_i(v_i)$ and can be controlled by integrating $\partial_t$ by parts. $J_1^*$ vanishes thanks to $\nabla^\varphi\cdot v=0$. We expand $J_2^*$
\begin{align*}
	J_2^* =  4\int_\Omega \partial_t^4 q(\partial_t \overline{\partial} \varphi \cdot  \partial_t^3 \partial_3 v) -4 \int_\Omega \partial_t^4 q \frac{\partial_\tau \varphi \partial_t\partial_3 \varphi}{\partial_3 \varphi}\partial_t^3\partial_3 v_\tau =:J_{21}^*+J_{22}^*.
\end{align*}
The boundary integral produced by integrating $\partial_3$ by parts in $J_{21}^*$ exactly cancels term $\int_\Sigma 4\partial_t^4 q \partial_t^3 \overline v\cdot  \partial_t \overline{\nabla}\psi$,
\begin{align*}
J_{21}^*-\int_\Sigma 4\partial_t^4 q \partial_t^3 \overline v\cdot  \partial_t \overline{\nabla}\psi= 4\int_\Omega \partial_3(\partial_t^4 q\partial_t \overline{\partial} \varphi) \cdot  \partial_t^3  v,
\end{align*}
which can be handled by integrating $\partial_t$ by parts under the time integral. The top order terms in $J_{22}^*+J_3^*$ vanish by carefully matching the terms thanks to the divergence-free condition
\begin{align*}
J_{22}^*+J_3^*&=4\int _\Omega 	 \partial_t^4 q \partial_t \partial_3 \varphi\partial_3^\varphi \partial_t^3  v_3 
+4 \int_\Omega \partial_t^4 q   \partial_t\partial_3 \varphi\big(\partial_\tau^\varphi \partial_t^3 v_\tau 
- \partial_\tau \partial_t^3 v_\tau\big)\\
&=4\int _\Omega 	 \partial_t^4 q \partial_t \partial_3 \varphi\partial_t ^3\big(\nabla^\varphi\cdot   v\big)  
-4\int _\Omega 	 \partial_t^4 q \partial_t \partial_3 \varphi\big(\big[\partial_t ^3,\nabla^\varphi\big]\cdot   v\big)
-4\int_\Omega \partial_t^4 q  \partial_t\partial_3 \varphi\partial_\tau \partial_t^3 v_\tau,
\end{align*}
the second term can be controlled by integrating $\partial_t$ by parts since $(\big[\partial_t ^3,\nabla^\varphi\big]\cdot   v\big)$ only contributes to lower order terms. The third term can be handled by integrating $\partial_t$ by parts and $\partial_\tau$ by parts to exchange the derivatives from $\partial_t^4q$ and $\partial_\tau\partial_t^3v_\tau$, where no boundary term produced due to our flat domain. Combining the analysis, we have 
\begin{align}
	\int_0^t\int_\Sigma \mathbf{Q}\mathcal{S}_1^* - \int_0^t\int_\Omega \partial_t^4 q \mathcal{C}_i(v_i)\partial_3 \varphi \leq \delta E_4^\kappa (t) +P(E_4^\kappa(0))+\int_0^tP(\sigma^{-1},E_4^\kappa(t)).
\end{align}

\section{Well-posedness of the nonlinear approximate system}
We have established the uniform-in$-\kappa $ estimate for our $\kappa-$approximate system $\ref{approxsys}$. In this section, we study the well-posedness of the approximate system for each fixed $\kappa>0$, which allows us to recover the well-posedness of the original system (\ref{sys}) by taking the limit $\kappa\to 0$. We use Picard iteration to construct a solution to the approximate system for each $\kappa>0$.
\subsection{Linearized approximate system}
Fix a $\kappa>0$, we start with $(v^{[0]},F_k^{[0]},\psi^{[0]})=(0,0,0)$ and $\psi^{[-1]}= \psi^{[0]}$. For any $n\in\mathbb{N}$, given $(v^{[n]},F_k^{[n]},q^{[n]},\psi^{[n]})$, we define $(v^{[n+1]},F_k^{[n+1]},q^{[n+1]},\psi^{[n+1]})$ to be the solution to the following linear system with coefficients only depending on $(v^{[n]},F_k^{[n]},q^{[n]},\psi^{[n]},\psi^{[n-1]})$
\begin{equation}\label{linearsys}
		\begin{cases}
			D_t^{\varphi^{[n]}} v^{[n+1]} +\nabla^{\varphi^{[n]}} q^{[n+1]}= (\mathfrak F_k^{[n]}\cdot \nabla^{\varphi^{[n]}}) F_k^{[n+1]} \qquad & \text{in } \Omega,\\
			\nabla^{\varphi^{[n]}} \cdot v^{[n+1]}=0 & \text{in } \Omega,\\
			D_t^{\varphi^{[n]}} F_j^{[n+1]}=(\mathfrak F_j^{[n]}\cdot \nabla^{\varphi^{[n]}}) v^{[n+1]}  & \text{in } \Omega,\\
			
			\partial_t \psi^{[n+1]}=v^{[n+1]}\cdot N^{[n]} &\text{on } \Sigma,\\
			q^{[n+1]}=-\sigma \overline{\nabla}\cdot\big(\frac{\overline{\nabla\psi^{[n]}}}{\sqrt{1+|\overline{\nabla}\psi^{[n]}|^2}}) + \kappa(1-\overline{\Delta})^2\psi ^{[n+1]}+ \kappa(1-\overline{\Delta})\partial_t\psi^{[n+1]} &\text{on } \Sigma,\\
			v_3^{[n+1]}=0 &\text{on } \Sigma_b,\\
			(v^{[n+1]},F_k^{[n+1]}, q^{[n+1]}, \psi^{[n+1]})|_{t=0} = (v_0, F_k^0, q_0, \psi_0).
		\end{cases}
	\end{equation}
	where $\mathfrak{F}^{[n]}_k$ is the modified deformation tensor defined by 
	\begin{align}
	\mathfrak{F}^{[n]}_{\tau k} &= F^{[n]}_{\tau k},	\qquad\tau=1,2,\\
	\mathfrak{F}^{[n]}_{3 k} &= F_{3k}^{[n]}+\mathfrak{R}_T(F_{1k}^{[n]}\partial_1\psi ^{[n]}+F_{2k}^{[n]}\partial_2\psi ^{[n]}-F_{3k}^{[n]})|_\Sigma .
	\end{align}
	and $\mathfrak{R}_T:L_T^2\big(H^s(\Sigma)\big)\to\mathfrak{R}_T:L_T^2\big(H^{s+1/2}(\Omega )\big) $ is a continuous extension operator such that $(\mathfrak R_T{g})|_\Sigma =g$ and 
	\begin{align} \label{extensioncontrol}
	||\mathfrak R_T{g}||_{s+\frac{1}{2}} \lesssim |g|_{s}.
	\end{align}
	The linearised derivatives are defined by 
	\begin{align}
		\partial_t^{\varphi^{[n]}}&=\partial_t-\frac{\partial_t {\varphi^{[n]}}}{\partial_3 {\varphi^{[n]}}}\partial_3,\\
		\partial_\tau^{\varphi^{[n]}}&=\partial_\tau-\frac{\partial_\tau {\varphi^{[n]}}}{\partial_3 {\varphi^{[n]}}}\partial_3, \qquad \tau=1,2,\\
		 \partial_3^{\varphi^{[n]}}&=\frac{1}{\partial_3 {\varphi^{[n]}}}\partial_3,
	\end{align}
	and the material derivative 
	\begin{align}
		D_t^{\varphi^{[n]}}
			=\partial_t+\overline{v}^{[n]}\cdot\overline{\partial}+\frac{1}{\partial_3 \varphi^{[n]}}(v^{[n]}\cdot\mathbf{N}^{[n-1]} -\partial_t{\varphi^{[n]}})\partial_3 ,
	\end{align}
	where $N^{[n]}= (-\partial_1\psi^{[n]} ,-\partial_2\psi^{[n]},1),\mathbf N^{[n]}= (-\partial_1\varphi^{[n]} ,-\partial_2\varphi^{[n]},1)$ with $\varphi^{[n]} = x_3+\chi (x_3)\psi^{[n]}.$ 
		\begin{remark}[Consistency of the material derivative and the kinematic boundary condition]	The defined linearised material derivative $D_t^{\varphi^{[n]}}$ no longer equals to $\partial_t^{\varphi^{[n]}}+v^{[n]}\cdot \nabla^{\varphi^{[n]}}$ due to the weight we set to $\partial_3$ to guarantee its consistency with the kinematic boundary condition. More precisely, we have $\frac{1}{\partial_3 \varphi^{[n]}}(v^{[n]}\cdot\mathbf{N}^{[n-1]} -\partial_t{\varphi^{[n]}})\partial_3$ vanished on the boundary.
		\end{remark}

	\begin{remark}[Propagation of the initial constraints]
		The initial constraint $\nabla^{\varphi^{[n]}}\cdot F_k^{[n]}$ and $F_{k}^{[n]}\cdot N^{[n]}= $ may not propagate due to the mismatch of $F_k^{[n]}$ and $F_k^{[n+1]}$ in the linearised system. We can still close the estimate without $\nabla^{\varphi^{[n]}}\cdot F_k^{[n]}$, however, the boundary constraint is crucial for our interior structure, otherwise, integrating $F_k^{[n]}\cdot \nabla^{\varphi^{[n]}}$ by parts would produce a uncontrollable boundary integral. We keep the constraint by defining the modified deformation tensor $\mathfrak{F}_k$ so that we have 
		\begin{align}
			\mathfrak{F}^{[n]}_k\cdot N^{[n]}=0 \qquad \text{on } \Sigma.
		\end{align}

	\end{remark}
	
	\begin{remark}[Recovering the approximate system]\label{recover}
		Due to the modified deformation tensor $\mathfrak{F}_k$, the linearised system does not directly recover the approximate system after taking the limit $n\to \infty$, but the following system
		\begin{equation}\label{limitsys}
		\begin{cases}
			D_t^{\varphi^{[\infty]}} v^{[\infty]} +\nabla^{\varphi^{[\infty]}} q^{[\infty]}= (\mathfrak F_k^{[\infty]}\cdot \nabla^{\varphi^{[\infty]}}) F_k^{[\infty]} \qquad & \text{in } \Omega,\\
			\nabla^{\varphi^{[\infty]}} \cdot v^{[\infty]}=0 & \text{in } \Omega,\\
			D_t^{\varphi^{[\infty]}} F_j^{[\infty]}=(\mathfrak F_j^{[\infty]}\cdot \nabla^{\varphi^{[\infty]}}) v^{[\infty]}  & \text{in } \Omega,\\
			
			\partial_t \psi^{[\infty]}=v^{[\infty]}\cdot N^{[\infty]} &\text{on } \Sigma,\\
			q^{[\infty]}=-\sigma \overline{\nabla}\cdot\big(\frac{\overline{\nabla\psi^{[\infty]}}}{\sqrt{1+|\overline{\nabla}\psi^{[\infty]}|^2}}) + \kappa(1-\overline{\Delta})^2\psi ^{[\infty]}+ \kappa(1-\overline{\Delta})\partial_t\psi^{[\infty]} &\text{on } \Sigma,\\
			v_3^{[\infty]}=0 &\text{on } \Sigma_b,\\
			(v^{[\infty]},F_k^{[\infty]}, q^{[\infty]}, \psi^{[\infty]})|_{t=0} = (v_0, F_k^0, q_0, \psi_0).
		\end{cases}
	\end{equation}
	We should show $\mathfrak{F}^{[\infty]}_{3k} = F_{3k}^{[\infty]}$, that is 
	\begin{align}
	\mathfrak{R}_T(F_{1k}^{[\infty]}\partial_1\psi ^{[\infty]}+F_{2k}^{[\infty]}\partial_2\psi ^{[\infty]}-F_{3k}^{[\infty]})|_\Sigma	 = 0.
	\end{align}
	It suffices to show $F_{1k}^{[\infty]}\partial_1\psi ^{[\infty]}+F_{2k}^{[\infty]}\partial_2\psi ^{[\infty]}-F_{3k}^{[\infty]}= 0 $ on $\Sigma$, then the above equality follows from the control (\ref{extensioncontrol}). Consider $D_t^{\varphi^{[\infty]}} F_j^{[\infty]}=(\mathfrak F_j^{[\infty]}\cdot \nabla^{\varphi^{[\infty]}}) v^{[\infty]}$ on $\Sigma$  and taking $\cdot N^{[\infty]}$,
	\begin{align*}
	D_t^{\varphi^{[\infty]}}(F_j^{[\infty]}\cdot N^{[\infty]}) 
	&= F_j^{[\infty]}\cdot D_t^{\varphi^{[\infty]}}N^{[\infty]} +(\mathfrak F_j^{[\infty]}\cdot \nabla^{\varphi^{[\infty]}}) v^{[\infty]}\cdot N^{[\infty]}\\
	&=F_j^{[\infty]}\cdot (\partial_t +\overline{v}^{[\infty]}\cdot \parbar) N^{[\infty]} + \big(\overline{\mathfrak F_j}^{[\infty]} \cdot \overline{\partial}+(\mathfrak{F}_j^{[\infty ]}\cdot N^{[\infty]})\partial_3\big) v^{[\infty]} \cdot N^{[\infty]}\\
	&=-(\overline{F_j}^{[\infty]}\cdot\parbar )(v^{[\infty]}\cdot N^{[\infty]})+(\overline{v}^{[\infty]}\cdot \parbar) N^{[\infty]}\cdot F_j^{[\infty]}+(\overline{F_j}^{[\infty]}\cdot\parbar )v^{[\infty]}\cdot N^{[\infty]}\\
	&=0.
	\end{align*}
	It follows that $F_{1k}^{[\infty]}\partial_1\psi ^{[\infty]}+F_{2k}^{[\infty]}\partial_2\psi ^{[\infty]}-F_{3k}^{[\infty]}= F_j^{[\infty]}\cdot N^{[\infty]}=0$ propagates from  the initial constraint $ F_j^{[\infty]}\cdot N^{[\infty]}|_{t=0}=0$.

	\end{remark}
	
	For simplicity of notations, for $n\in \mathbf{N}$, we define
	\begin{align*}
		(v,F_k,q,\psi)&=(v^{[n+1]},F_k^{[n+1]},q^{[n+1]},\psi^{[n+1]}),\\
		(\mathring v,\mathring F_k,\mathring {\mathfrak F_k}, \mathring q,\mathring \psi)&=(v^{[n]},F_k^{[n]},\mathfrak F_k^{[n]},q^{[n]},\psi^{[n]}),\\
		\dot \psi &= \psi^{[n-1]}.
	\end{align*}
	then the system (\ref{linearsys}) becomes
	\begin{equation}\label{linearsys2}
		\begin{cases}
			D_t^{\mathring \varphi} v +\nabla^{\mathring \varphi} q= (\mathring {\mathfrak F_k}\cdot \nabla^{\mathring \varphi}) F_k \qquad & \text{in } \Omega,\\
			\nabla^{\mathring \varphi} \cdot v=0 & \text{in } \Omega,\\
			D_t^{\mathring \varphi} F_j=(\mathring {\mathfrak F_j}\cdot \nabla^{\mathring \varphi}) v  & \text{in } \Omega,\\
			
			\partial_t \psi=v\cdot \mathring N &\text{on } \Sigma,\\
			q=-\sigma \overline{\nabla}\cdot\big(\frac{\overline{\nabla\mathring \psi}}{\sqrt{1+|\overline{\nabla}\mathring \psi|^2}}) + \kappa(1-\overline{\Delta})^2\psi+ \kappa(1-\overline{\Delta})\partial_t\psi &\text{on } \Sigma,\\
			v_3=0 &\text{on } \Sigma_b,\\
			(v,F_k, q, \psi)|_{t=0} = (v_0, F_k^0, q_0, \psi_0).
		\end{cases}
	\end{equation}
Now, we assume the bounds for the basic state $(\mathring v,\mathring F_k,\mathring {\mathfrak F_k}, \mathring q,\mathring \psi)$ and $
		\dot \psi$: there exists $\mathring K_0>0, T_\kappa>0$ such that, for $t\in [0,T_\kappa]$,
	\begin{align}
	 \sum_{l=0}^4||\partial_t^l( \mathring v, \mathring F_k)||^2_{4-l}+ \sum_{l=0}^4|\sqrt\kappa \partial_t^l \mathring \psi|^2_{6-l}+\int_0^t|\sqrt\kappa \partial_t^5 \mathring\psi|^2_1\leq \mathring K_0	.
	\end{align}
	It follows that we have control for the modified deformation tensor $\mathring {\mathfrak F_k}$,
	\begin{align}	
	\int_0^t||\partial_t^l \mathring {\mathfrak F_k}||_{4-l}^2 \leq C(\mathring K_0).
	\end{align}
	
	Before going into the analysis, we introduce the linearised transport theorem, in which we have additional mismatched terms due to our definition $D_t^{\mathring\varphi}$.
	\begin{lemma}[Linearised transport Theorem]\label{lintransp}
		 Let $f=f(t,x), x\in \Omega.$
		\begin{align}
		\frac{1}{2}\frac{d}{dt}\int_\Omega |f|^2\partial_3 \mathring\varphi =\int_\Omega (D_t^{\mathring\varphi} f)f\partial_3 \mathring\varphi +\frac{1}{2}\int_\Omega|f|^2(\nabla^{\mathring\varphi}\cdot\mathring v)\partial_3\mathring \varphi+\frac{1}{2}\int_\Omega |f|^2\bigg(\partial_3(\overline{\mathring v}\cdot \overline{\nabla})(\mathring\varphi-\dot\varphi)\bigg).
		\end{align}
	\end{lemma}
	The proof is straightforward computation by expanding $D_t^{\mathring\varphi}$, where standard integration by parts gives us the $\frac{1}{2}$ weighted terms by observing the symmetry. We refer to \cite{luo2022compressible} for details. 
		\subsection{Well-posedness of the linearised approximate system}

	We apply the classical Galerkin’s method to construct the weak solution in $L^2([0,T]\times \Omega)$ to the linearized problem (\ref{linearsys2}), which is actually a strong solution by the argument in \cite{metivier2012small}. 	
	Let $\{\mathbf{e}_j\}_{j=1}^\infty\subseteq C^\infty (\Omega)$ be an orthogonal basis of $L^2(\Omega)$ and $H_0^1(\Omega).$ We define $ U=(q, v , \mathbf {F}_k) $, then the linearised system (\ref{linearsys2}) can be expressed in terms of $\mathring U$ by 
	\begin{align}\label{matrixform}
	A_0(\mathring U)\partial_t U+\sum_{i=1}^3 A_i(\mathring U)\partial_i U=0,
	\end{align}
	where, for $\tau = 1 ,2,$
	\begin{align*}
		A_0(\mathring U) &= 
		\begin{bmatrix}
 		0 &  \overrightarrow0^\intercal & \overrightarrow0^\intercal\\
 		\overrightarrow0 & \mathbf I_3 & \mathbf O_3\\
 		\overrightarrow0 & \mathbf O_3 & \mathbf I_3
		 \end{bmatrix} , 
		 A_\tau (\mathring U) = 
		 \begin{bmatrix}
 		0 &  \overrightarrow e^\intercal_\tau & \overrightarrow0^\intercal\\
 		\overrightarrow e_\tau & \mathring v_\tau \mathbf I_3 & -\mathring {\mathfrak{F}}_{\tau k}\mathbf I_3\\
 		\overrightarrow0 & -\mathring{\mathfrak{F}}_{\tau k}\mathbf I_3 & \mathring v_\tau \mathbf I_3 
		 \end{bmatrix},\\
		 A_3(\mathring U) &= \frac{1}{\partial _3 \mathring \varphi}
		 \begin{bmatrix}
		 	0 & \mathring{ \mathbf N }^\intercal & \overrightarrow0^\intercal\\
		 	\mathring {\mathbf N} & (\mathring v \cdot \dot{\mathbf N}-\partial_t \mathring \varphi)\mathbf I_3 & -(\mathfrak F_k\cdot  \mathring {\mathbf N})\mathbf I_3\\
		 	\overrightarrow 0 & -(\mathfrak F_k\cdot  \mathring {\mathbf N})\mathbf I_3 & (\mathring v \cdot \dot{\mathbf N}-\partial_t \mathring \varphi)\mathbf I_3
		 \end{bmatrix}.
	\end{align*}
	For $m\in \mathbb N$, we define the Galerkin sequence $(U^m(t,x), \psi^m(t,\xbar)) $ by
	\begin{align}
	U_j^m(t,x) &= \sum_{j=1}^m U^m_{lj}(t)\mathbf{e}_l(x), \\
	\psi^m(t,\xbar) &= \sum_{l=1}^m\psi _l^m(t)\mathbf{e}_l(\xbar,0),
	\end{align}
	with the boundary conditions
	\begin{align}
		\partial_t\psi^m &= - U_2^m\partial_1\mathring \psi-U_3^m\partial_2\mathring \psi+U_4^m,\\
		U_1^m &= -\sigma \overline{\nabla}\cdot\big(\frac{\overline{\nabla\mathring \psi}}{\sqrt{1+|\overline{\nabla}\mathring \psi|^2}}) + \kappa(1-\overline{\Delta})^2\psi^m+ \kappa(1-\overline{\Delta})\partial_t\psi^m.
	\end{align}
	Let
	\begin{align}
	\phi_i(t,x) = \sum_{l=1}^m\phi_{il}(t)\mathbf{e}_l(x).
	\end{align}
	We test (\ref{matrixform})  with $\phi_i$ and integrating $\partial_i$ by parts to get
	\begin{align*}
		\int_\Omega A^{ij}_0(\mathring U)\partial_t U_j\phi_i- \sum_{s=1}^3 \int_\Omega \partial_s (A^{ij}_s(\mathring U) \phi_i)U_j +\int_\Sigma A^{ij}_3(\mathring U) U_j\phi_i=0.
	\end{align*}
	Plugging in the Galerkin sequence and taking $\phi_i= \mathbf{e}_i$, we have an ODE of $U_{lj}^m(t)$
	\begin{align*}
		&\bigg(\int_\Omega A^{ij}_0(\mathring U) \mathbf{e}_l\mathbf{e}_i\bigg)(U^m_{lj})'(t) - \bigg(\sum_{s=1}^3 \int_\Omega \partial_s (A^{ij}_s(\mathring U) \mathbf{e}_i)\mathbf{e}_{l}\bigg) U_{lj}^m(t)\\
		&= -\int_\Sigma\partial_t\psi^m \mathbf{e}_1 -\int_\Sigma U^m_1(-\partial_1\mathring \psi \mathbf{e}_2-\partial_2\mathring \psi \mathbf{e}_3+ \mathbf{e}_4)\\
		&=\sigma \int_\Sigma\big(\frac{\overline{\nabla\mathring \psi}}{\sqrt{1+|\overline{\nabla}\mathring \psi|^2}})\cdot \overline{\nabla}(\phi_\nu\cdot \mathring N)-\kappa \int_\Sigma(1-\overline{\Delta})\psi^m(1-\overline{\Delta})(\phi_\nu\cdot \mathring N)-\kappa \int_\Sigma \partial_t\psi^m(\phi_\nu\cdot \mathring N)\\
		&\quad-\int_\Sigma\partial_t\psi^m \mathbf{e}_1+\kappa \int_\Sigma \overline{\nabla}\partial_t\psi^m\overline{\nabla}(\phi_\nu\cdot \mathring N),
	\end{align*}
	where $\phi_\nu = -\partial_1\mathring \psi \mathbf{e}_2-\partial_2\mathring \psi \mathbf{e}_3+ \mathbf{e}_4$, 
	with the initial data
	\begin{align*}
		U_{lj}^m(0) = \int_\Omega U_j^m(0,x) \mathbf{e}_l(x).
	\end{align*}
	The local well-posedness and existence of the above system are guaranteed by standard ODE theories.
	
	Now we derive the uniform-in-$m$ estimate of the Galerkin sequence  $(U^m(t,x), \psi^m(t,\xbar)) $, which guarantees a weak solution to (\ref{matrixform}) by taking the limit $m\to \infty$ subject to a subsequence. Testing (\ref{matrixform}) with $U^m$ and integrating $\partial_i$ by parts, 
	\begin{align*}
	\frac{1}{2}\frac{d}{dt}\int _\Omega A_0(\mathring U)U^m(U^m)^\intercal = \frac{1}{2}\int_\Omega \partial_tA_0(\mathring U)U^m(U^m)^\intercal+\frac{1}{2}\sum_{i=1}^3 \int_\Omega \partial_iA_i(\mathring U) U^m(U^m)^\intercal-\int_\Sigma A_3(\mathring U) U^m(U^m)^\intercal,
	\end{align*}
	where the $\frac{1}{2}$ weight comes from the symmetry of the coefficient matrix.  The interior integral can be directly controlled by $C(\mathring K_0)||U^m||_0^2$ since there is no derivative landed on $U^m$. Invoking the boundary conditions for the Galerkin sequence, the boundary integral gives
	\begin{align*}
		-\int_\Sigma A_3(\mathring U) U^m(U^m)^\intercal  
		&= -2\int_\Sigma \bigg(-\sigma \overline{\nabla}\cdot\big(\frac{\overline{\nabla\mathring \psi}}{\sqrt{1+|\overline{\nabla}\mathring \psi|^2}}) + \kappa(1-\overline{\Delta})^2\psi^m+ \kappa(1-\overline{\Delta})\partial_t\psi^m\bigg)\partial_t\psi^m,
	\end{align*}
	where the last two terms contribute to energy terms
	\begin{align*}
	-2 \int_\Sigma 	\kappa(1-\overline{\Delta})^2\psi^m\partial_t\psi^m 
	&= -\frac{d}{dt}\int_\Sigma|\sqrt\kappa(1-\overline{\Delta})\psi^m|^2,\\
	-2\int_\Sigma \kappa(1-\overline{\Delta})\partial_t\psi^m\partial_t\psi^m
	&=-\frac{d}{dt}\int_0^t\int_\Sigma |\sqrt\kappa\japbr\partial_t\psi^m|^2,
	\end{align*}
	and the first term can be directly controlled,
	\begin{align*}
	2\int_0^t\int_\Sigma \sigma \overline{\nabla}\cdot\big(\frac{\overline{\nabla\mathring \psi}}{\sqrt{1+|\overline{\nabla}\mathring \psi|^2}})	\partial_t\psi^m\leq \delta\int_0^t \int_\Sigma |\sqrt\kappa \partial_t\psi^m|_1^2+C(\mathring K_0)t.
	\end{align*}
	Summarising the analysis, we obtained the uniform-in-$n$ energy estimate 
	\begin{align}
		E^m(t)\leq E^m(0)+\int_0^t C(\mathring K_0,\kappa^{-1})E^m(s)ds,
	\end{align}
	where the energy is defined by
	\begin{align}
		E^m(t) = ||U^m||_0^2 +|\sqrt\kappa \psi^m|_2^2+\int_0^t |\sqrt{\kappa}\partial_t\psi|_1^2.
	\end{align}
	Hence, by Gronwall's inequality, there exists $T>0$ independent of $m$ such that, for $t\in[0,T]$ 
	\begin{align}
		E^m(t)\leq C(\mathring K_0,\kappa^{-1})E(0).
	\end{align}

	By Eberlein-Sˇmulian theorem, the Galerkin sequence weakly converges subject to a subsequence,
	\begin{align}
		U^{m_k}\rightharpoonup(q,v,\mathbf F_k)& \text{ in }L^2([0,T];L^2(\Omega)),\\
		\psi^{m_k}\rightharpoonup\psi&  \text{ in }L^2([0,T];H^2(\Omega)),\\
		\partial_t\psi^{m_k}\rightharpoonup\partial_t\psi&  \text{ in }L^2([0,T];H^1(\Omega)).
	\end{align}
	Thus, we have proved the existence of a weak solution to the linearised system (\ref{linearsys2}) and the uniqueness directly follows from the energy estimate of $E^m(t).$

	\subsection{Higher-order estimates of the linearized approximate system}\label{highlin}
	Now we prove the higher-order energy estimates of the linearised system (\ref{linearsys2}), with which we can proceed with the Picard iteration. 
	
	\begin{prop} \label{higherorderlinearised}Fix $\kappa>0$, there exists some $T_\kappa>0$ such that, for $t\in [0,T_\kappa],$
	\begin{align}
		\mathring{E}_\kappa^4(t)\leq C(\kappa^{-1},\mathring K_0)\mathring{E}_\kappa^4(0),
	\end{align}
	where the energy $\mathring{E}_\kappa^4(t)$ is defined by 
	\begin{align}
		\mathring{E}_\kappa^4(t):= \sum_{l=0}^4||\partial_t^l(  v, F_k)||^2_{4-l}+ \sum_{l=0}^4|\sqrt\kappa \partial_t^l  \psi|^2_{6-l}+\int_0^t|\sqrt\kappa \partial_t^5 \psi|^2_1.
	\end{align}
	\end{prop}
	We no longer expect energy terms from the surface tension $\sigma \mathcal{H}(\mathring\psi)$, since it was set to be fully given. However, the artificial viscosity terms $\kappa(1-\overline{\Delta})^2\psi$ and $\kappa(1-\overline{\Delta})\partial_t\psi$ can give us higher regularity as we are solving the equation for fixed $\kappa$ and no longer need to keep the estimate uniform-in-$\kappa$.  Moreover, for the solved $\mathring \psi$, we can have a bound for $\sum_{l=0}^3|\partial_t^l \psi|_{7-l}$ by observing the biharmonic equation on the boundary
		\begin{align*}
		\overline{\Delta}^2\mathring \psi
		&=(1-\overline{\Delta})^2\mathring \psi-\mathring\psi+2\overline{\Delta}\mathring\psi\\
		&=\kappa^{-1} \bigg(-\sigma \overline{\nabla}\cdot\big(\frac{\overline{\nabla\dot \psi}}{\sqrt{1+|\overline{\nabla}\dot \psi|^2}}) + \kappa(1-\overline{\Delta})\partial_t\mathring \psi-\mathring q\bigg)-\mathring\psi+2\overline{\Delta}\mathring\psi,
	\end{align*}
	where we can control the right-hand side up to 3 derivatives. 

 \subsubsection{Div-curl estimate}
 By lemma \ref{hodge}, for $l=0,1,2,3,$ we have 
 \begin{align}
		||\partial_t^l v||_{4-l}^2&\lesssim C(\mathring K_0)\big(||\nabla^{\mathring\varphi}\times \partial_t^l v||_{3-l}^2+||\overline{\partial}^{4-l} \partial_t^l v||_0^2+|| \partial_t^lv||_0^2\big),\\
		||\partial_t^lF_j||_{4-l}^2&\lesssim C(\mathring K_0)\big(||\nabla^{\mathring\varphi}\cdot \partial_t^l F_j||_{3-l}^2+||\nabla^{\mathring\varphi}\times \partial_t^l F_j||_{3-l}^2+||\overline{\partial}^{4-l} \partial_t^l F_j||_0^2+||\partial_t^l F_j||_0^2\big).
	\end{align}
	
	We have studied the $L^2$ part in the uniform-in-m estimate of the Glerkin sequence. For the divergence part, we still have $\nabla^{\mathring\varphi}\cdot v=0$, but we no longer have the divergence-free constraint for $F_k$. Taking $\nabla^{\mathring\varphi}$ to the linearised evolution equation of $F_j$, we get
	\begin{align}\label{divF}
	D_t^{\mathring \varphi} (\nabla^{\mathring\varphi}\cdot F_j)=(\partial^{\mathring\varphi}_i\mathring {\mathfrak F_j}\cdot \nabla^{\mathring \varphi}) v_i+[D_t^{\mathring\varphi},\nabla^{\mathring\varphi}\cdot ]F_j,
	\end{align}
	where $[D_t^{\mathring\varphi},\nabla^{\mathring\varphi}\cdot ]F_j = -\partial_i^{\mathring\varphi}\mathring v\cdot \nabla^{\mathring\varphi}F_{ij}-(\partial_i^{\mathring\varphi}\partial_t(\mathring\varphi-\dot\varphi))\partial_3 ^{\mathring\varphi}F_{ij}$. Note that there is only one derivative lands on $v$ and $F_j$ and the energy terms of $\dot\psi,\mathring \psi$ can be controlled by $\mathring K$. Hence, we can obtain the control of $||\nabla^{\mathring\varphi}\cdot \partial_t^lF_j||_{3-l}$ by testing the $\partial_t^l\partial^{3-l}$-differentiated (\ref{divF}) with $\partial^{3-l}(\nabla^{\mathring\varphi}\cdot \partial_t^lF_j)$
	\begin{align}
		\frac{1}{2}\frac{d}{dt}||\nabla^{\mathring\varphi}\cdot \partial_t^lF_j||_{3-l}^2\lesssim C(\mathring K_0,\kappa^{-1})(||\partial_t^l v||_{4-l},||\partial_t^l v||_{4-l}).
	\end{align}
	For the vorticity part, we take the curl operator $\nabla^{\mathring\varphi}\times $ and $\partial^3$ to the evolution equations of $v,F_j$ to get 
	\begin{align}
		D_t^{\mathring \varphi}(\partial^{3}\nabla^{\mathring\varphi}\times v)&= (\mathring {\mathfrak F}_k\cdot\nabla^{\mathring \varphi})(\partial^3\nabla^{\mathring\varphi}\times F_k)+ R_1,\\
		D_t^{\mathring \varphi}(\partial^{3}\nabla^{\mathring\varphi}\times F_j)&= (\mathring{\mathfrak F}_j\cdot\nabla^{\mathring \varphi})(\partial^3\nabla^{\mathring\varphi}\times v)+ R_2,
	\end{align}
	where $R_1= \partial^3[\nabla^{\mathring\varphi}\times,D_t^{\mathring\varphi}]v+\partial^3 [\nabla^{\mathring\varphi}\times,(\mathring {\mathfrak F}_k\cdot\nabla^{\mathring \varphi})]F_k+[\partial^3,D_t^{\mathring\varphi}]\nabla^{\mathring\varphi}\times v+[\partial^3,(\mathring {\mathfrak F}_k\cdot\nabla^{\mathring \varphi})] (\nabla^{\mathring\varphi}\times F_k)$ and $R_2 = \partial^3[\nabla^{\mathring\varphi}\times,D_t^{\mathring\varphi}]F_j+\partial^3 [\nabla^{\mathring\varphi}\times,(\mathring {\mathfrak F}_j\cdot\nabla^{\mathring \varphi})]v+[\partial^3,D_t^{\mathring\varphi}]\nabla^{\mathring\varphi}\times F_j+[\partial^3,(\mathring {\mathfrak F}_j\cdot\nabla^{\mathring \varphi})] (\nabla^{\mathring\varphi}\times v)$.
	The top order terms $(\mathring {\mathfrak F}_k\cdot\nabla^{\mathring \varphi})(\partial^3\nabla^{\mathring\varphi}\times F_k)$ and $(\mathring{\mathfrak F}_j\cdot\nabla^{\mathring \varphi})(\partial^3\nabla^{\mathring\varphi}\times v)$ can be cancelled due to the symmetric structure when we test the first equation with  $(\partial^{3}\nabla^{\mathring\varphi}\times v)$ and the second equation with $(\partial^{3}\nabla^{\mathring\varphi}\times F_j)$. It follows that
	\begin{align}
		\frac{1}{2}\frac{d}{dt}\big(||\partial^3\nabla^{\mathring\varphi}\times v||_0^2+||\partial^{3}\nabla^{\mathring\varphi}\times F_k||_0^2 \big)\lesssim C(\mathring K_0,\kappa^{-1})(||F_k||_{4},||v||_{4}).
	\end{align}
	We can apply parallel arguments to treat the case with $\partial_t^l\partial^{3-l}$ to get
		\begin{align}
		\frac{1}{2}\frac{d}{dt}\big(||\nabla^{\mathring\varphi}\times \partial_t^lv||_{3-l}^2+||\nabla^{\mathring\varphi}\times \partial_t^lF_k||_{3-l}^2 \big)\lesssim C(\mathring K_0,\kappa^{-1})(||\partial_t^l F_k||_{4-l},||\partial_t^l v||_{4-l}).
	\end{align}

	\subsubsection{Tangential Estimates}\label{lintan}
	The tangential analysis is parallel to the uniform energy estimate of the nonlinear approximate system. As we mentioned, it is even easier as we can freely use the $\sqrt\kappa-$weighted energy to control the terms without worrying about the attachment of $\kappa$. Let $D^\alpha = \partial_t^{\alpha_0}\partial_1^{\alpha_1}\partial_2^{\alpha_2}$ with $|\alpha|\leq 4$. We introduce Alinhac's good unknowns for the $D^\alpha$ differentiated (\ref{linearsys2}), for which we should expect no difference to the ones for the nonlinear system except for the weights of $\partial_3 $ in $D_t^{\mathring \varphi}$,
	\begin{align}
		D^\alpha \partial_i^{\mathring \varphi} f
		=\partial^\varphi_i(D^\alpha f-D^\alpha\mathring \varphi\partial_3^{\mathring \varphi} f) +\mathring{\mathcal{C}}_i(f),
	\end{align}
	where, for $\tau =1,2, |\beta|=1$,
	\begin{align}
		\mathring{\mathcal{C}}_\tau (f)&= D^\alpha \mathring \varphi\partial_\tau^{\mathring \varphi}\partial_3^{\mathring \varphi} f -[D^\alpha,\frac{\partial_\tau \mathring \varphi}{\partial_3 \mathring \varphi },\partial_3 f]-\partial_3 f[D^\alpha,\partial_\tau \mathring\varphi,\frac{1}{\partial_3 \mathring\varphi}]+\partial_3 f \partial _\tau \mathring \varphi[D^{\alpha-\beta},\frac{1}{(\partial_3 \mathring \varphi)^2}]D^\beta \partial_3 \mathring\varphi,\\
		\mathring{\mathcal{C}}_3(f)&=D^\alpha \mathring \varphi (\partial_3 ^{\mathring \varphi})^2 f +[D^\alpha,\frac{1}{\partial_3 \mathring \varphi},\partial_3 f]- \partial_3 f [D^{\alpha-\beta},\frac{1}{(\partial_3 \mathring \varphi)^2}]D^\beta \partial_3 \mathring\varphi,
	\end{align}
	and 
	\begin{align}
		D^\alpha D_t^{\mathring \varphi} f
		=D_t^{\mathring \varphi} (D^\alpha f- D^\alpha \mathring \varphi\partial_3 ^{\mathring \varphi} f)+\mathring{\mathcal{D}}(f),
	\end{align}
	where 
	\begin{align}
			\mathring{\mathcal{D}}(f)&=D^\alpha \mathring \varphi D_t^{\mathring\varphi} \partial_3^{\mathring \varphi} f +[D^\alpha, \overline{\mathring v}]\cdot \overline{\nabla}f+[D^\alpha, \frac{1}{\partial_3 \mathring\varphi}(\mathring v\cdot \dot{\mathbf N}-\partial_t \mathring\varphi),\partial_3 f]+[D^\alpha, \frac{1}{\partial_3 \mathring \varphi},\mathring v\cdot \dot {\mathbf{N}}-\partial_t \mathring \varphi]\partial_3 f\nonumber\\
			&\quad-(\mathring v\cdot \dot{\mathbf{N}} -\partial_t \mathring \varphi)\partial_3 f[D^{\alpha-\beta},\frac{1}{(\partial_3 \mathring \varphi)^2}]D^\beta \partial_3 \mathring \varphi + \frac{1}{\partial_3 \mathring \varphi}\partial_3 f[D^\alpha,\mathring v]\dot{\mathbf{N}}.
	\end{align}
	
	We define 
	\begin{align}
		\mathring{\mathbf V}=  D^\alpha v-D^\alpha\mathring \varphi\partial_3^{\mathring \varphi} v,
		\qquad \mathring{\mathbf F_j}=D^\alpha F_j-D^\alpha\mathring \varphi\partial_3^{\mathring \varphi} F_j,
		\qquad \mathring{\mathbf Q}=  D^\alpha q-D^\alpha\mathring \varphi\partial_3^{\mathring \varphi} q,
	\end{align}
	then the $D^\alpha$-differentiated (\ref{linearsys2}) can be written in terms of the AGUs as follows

\begin{align}\label{linAGUsys}
\begin{cases}
	D_t^{\mathring \varphi} \mathring{\mathbf{V}}_i + \partial_i^{\mathring \varphi}  \mathring{\mathbf{Q}}= (\mathring {\mathfrak F}_k\cdot \nabla^{\mathring \varphi} ) \mathring{\mathbf{F}}_{ik}+\mathring{\mathcal{R}}_i^1 \qquad & \text{in } \Omega,\\
	\partial_k^{\mathring \varphi} \mathring{ \mathbf{V}}_k=-\mathring{\mathcal{C}}_k(v_k) & \text{in } \Omega,\\
	D_t^{\mathring \varphi}\mathring { \mathbf{F}}_{ij}=(\mathring{\mathfrak F}_j\cdot \nabla^{\mathring \varphi})\mathring{\mathbf{V}}_i+\mathring{\mathcal  R}_{ij}^2  & \text{in } \Omega, \\
	\mathring{\mathbf{Q}} = -\sigma D^\alpha\big(\overline{\nabla}\cdot \frac{\overline{\nabla}\mathring \psi}{|\mathring N|}\big)
	+ \kappa D^\alpha(1-\overline{\Delta})^2\psi + \kappa D^\alpha(1-\overline{\Delta})\partial_t\psi 
	-D^\alpha \mathring \psi \partial_3 q   &\text{on }\Sigma , \\
	\partial_t D^\alpha \psi =\mathring{\mathbf{V}}\cdot \mathring N-(\overline{v}\cdot \overline{\nabla})  D^\alpha \mathring \psi+\mathring{\mathcal{S}}_1  &\text{on }\Sigma. 
\end{cases} 
\end{align}
where 
\begin{align}
\mathring{\mathcal{R}}_i^1&=\mathring{\mathfrak F}_{lk}\mathring{\mathcal{C}}_l(F_{ik})+[D^\alpha, \mathring{\mathfrak F}_{lk}] \partial_l^{\mathring\varphi} F_{ik}-\mathring{\mathcal{D}}(v_i)-\mathring{\mathcal{C}}_i(q), \\
\mathring{\mathcal{R}}_{ij}^2&=\mathring {\mathfrak F}_{kj}\mathring {\mathcal{C}}_k(v_i)+[D^\alpha,\mathring{\mathfrak F}_{kj}]\partial_k^{\mathring\varphi} v_i-\mathring{\mathcal{D}}(F_{ij}),\\
\mathring{\mathcal{S}}_1&=D^\alpha \mathring \psi \partial _3 v\cdot \mathring N -[D^\alpha,\overline{v}\cdot,\overline{\partial}\mathring\psi].
\end{align}
\textbf{Part 1: Interior structure}

In the analysis of linearised equations, there are additional mismatched terms produced in the linearised Reynold's transport theorem and interior integrals attached with $\nabla^{\mathring \varphi}\cdot \mathring{\mathfrak F}_j$ when integrating $ \mathring{\mathfrak F}_j\cdot\nabla^{\mathring \varphi}$ by parts, but they are lower order terms and does not cause extra difficulties. Testing the first equation of (\ref{linAGUsys}) with $\mathring{\mathbf{V}}_i\partial_3\mathring\varphi$, we get 
\begin{align}\label{linint}
	\frac{1}{2}\frac{d}{dt}\int_\Omega |\mathring{\mathbf{V}}_i|^2\partial_3\mathring\varphi 
	&=-\int_\Omega \mathring{\mathbf{V}}_i\partial_i^{\mathring \varphi}  \mathring{\mathbf{Q}}\partial_3\mathring\varphi
	+\int_\Omega \mathring{\mathbf{V}}_i(\mathring{\mathfrak F}_k\cdot \nabla^{\mathring \varphi} ) \mathring{\mathbf{F}}_{ik}\partial_3\mathring\varphi
	+\int_\Omega \mathring{\mathbf{V}}_i\mathring{\mathcal{R}}_i^1\partial_3\mathring\varphi \nonumber\\
	&\quad+\frac{1}{2}\int_\Omega|\mathring{\mathbf{V}}_i|^2\bigg((\nabla^{\mathring\varphi}\cdot\mathring v)\partial_3\mathring \varphi+\partial_3(\overline{\mathring v}\cdot \overline{\nabla})(\mathring\varphi-\dot\varphi)\bigg),
\end{align}
where the last two terms can by directly controlled by $C(\mathring K_0)\mathring E_4^\kappa (t)$. The second term contributes to the energy of $\mathring {\mathbf {F}}$ by integrating $(\mathfrak F_k\cdot \nabla^{\mathring \varphi} )$ by parts and invoking the second equation of (\ref{linAGUsys}),
\begin{align*}
	&\int_\Omega \mathring{\mathbf{V}}_i(\mathring{\mathfrak F}_k\cdot \nabla^{\mathring \varphi} ) \mathring{\mathbf{F}}_{ik}\partial_3\mathring\varphi \\
	&=-\int_\Omega (\mathring{\mathfrak F}_k\cdot \nabla^{\mathring \varphi} )\mathring{\mathbf{V}}_i \mathring{\mathbf{F}}_{ik}\partial_3\mathring\varphi -\int_\Omega (\nabla^{\mathring \varphi}\cdot \mathring{\mathfrak F}_k)\mathring{\mathbf{V}}_i \mathring{\mathbf{F}}_{ik}\partial_3\mathring\varphi \\
	&=-\frac{1}{2}\frac{d}{dt}\sum_k\int_\Omega |\mathring { \mathbf{F}}_{ik}|^2 \partial_3\mathring\varphi 
	+\frac{1}{2}\sum_k\int_\Omega|\mathring{\mathbf{F}}_{ik}|^2\bigg((\nabla^{\mathring\varphi}\cdot\mathring v)\partial_3\mathring \varphi+\partial_3(\overline{\mathring v}\cdot \overline{\nabla})(\mathring\varphi-\dot\varphi)\bigg)
	+\int_\Omega \mathring{\mathcal  R}_{ij}^2\mathring{\mathbf{F}}_{ik}\partial_3\mathring\varphi\\
	&\quad-\int_\Omega (\nabla^{\mathring \varphi}\cdot \mathring{\mathfrak F}_k)\mathring{\mathbf{V}}_i \mathring{\mathbf{F}}_{ik}\partial_3\mathring\varphi,
\end{align*}
where the boundary integral vanishes due to $\mathring{\mathfrak{F}_k}\cdot \mathring N=0$ on the boundary and the last three terms can be directly controlled. The first term of (\ref{linint}) give rise to the boundary integral by integrating $\partial_i^{\mathring\varphi}$ by parts, 
\begin{align}
-\int_\Omega \mathring{\mathbf{V}}_i\partial_i^{\mathring \varphi}  \mathring{\mathbf{Q}}\partial_3\mathring\varphi	 = \int _\Omega \mathring {\mathcal C}_i(v_i)\mathring{\mathbf Q}\partial_3\mathring\varphi-\int_\Sigma \mathring{\mathbf Q}\big(\mathring {\mathbf{V}}\cdot \mathring N\big).
\end{align}
\textbf{Part 2: Boundary structure}\\
Most of the boundary integrals can be handled by arguments parallel to the ones in section \ref{prioritangential} or directly controlled as we have an even higher regularity of $\psi$. The only parts that we should pay attention to are the terms that were handled using symmetry since the linearisation breaks the symmetric structure. 
\begin{align}\label{linpart2}
	\int_\Sigma \mathring{\mathbf Q}\big(\mathring {\mathbf{V}}\cdot \mathring N\big) 
	=
	\int_\Sigma \mathring{\mathbf Q}\partial_t D^\alpha \psi 
	+\int_\Sigma \mathring{\mathbf Q}(\overline{v}\cdot \overline{\nabla})  D^\alpha \mathring \psi
	-\int_\Sigma \mathring{\mathbf Q}\mathring{\mathcal{S}}_1.
\end{align}
The first term gives
\begin{align*}
	\mathring{ST}+\mathring{RT} 
	=\int_\Sigma \big( -\sigma D^\alpha\big(\overline{\nabla}\cdot \frac{\overline{\nabla}\mathring \psi}{|\mathring N|}\big)
	+ \kappa D^\alpha(1-\overline{\Delta})^2\psi + \kappa D^\alpha(1-\overline{\Delta})\partial_t\psi \big)\partial_t D^\alpha \psi
	-\int_\Sigma D^\alpha \mathring \psi \partial_3 q \partial_t D^\alpha \psi,
\end{align*}
where the $\mathring{RT}$ can be directly controlled even when $D^\alpha=\partial_t^4$, since we can control $\partial_t^4\mathring \psi$ by $\mathring K_0$ and $\partial_t^5\psi$ by the $\kappa$-weighted energy under time integral. We expect the last two terms of $\mathring {ST}$ to contribute to the boundary regularity and we directly control the first term under time integral,
\begin{align*}
	\int_\Sigma\big( \kappa D^\alpha(1-\overline{\Delta})^2\psi + \kappa D^\alpha(1-\overline{\Delta})\partial_t\psi \big)\partial_t D^\alpha \psi = -\frac{1}{2}\frac{d}{dt}  \int_\Sigma |\sqrt\kappa D^\alpha(1-\overline\Delta)\psi|^2 -\frac{d}{dt}\int _0^t\int_\Sigma|\sqrt{\kappa}\partial_tD^\alpha \japbr\psi|^2,
	\end{align*}
	and 
	\begin{align*}
		-\int_0^t\int_\Sigma \sigma D^\alpha\big(\overline{\nabla}\cdot \frac{\overline{\nabla}\mathring \psi}{|\mathring N|}\big)\partial_t D^\alpha \psi \leq \frac{\delta}{\kappa}\int_0^t|\sqrt\kappa \partial_t D^\alpha \psi|^2+\frac{\sigma}{\kappa}\int_0^t C(\mathring {K}_0).
	\end{align*}
	The second term of (\ref{linpart2}) gives 
	\begin{align*}
		&\int_\Sigma \mathring{\mathbf Q}(\overline{v}\cdot \overline{\nabla})  D^\alpha \mathring \psi \\
		&= \int_\Sigma \big( -\sigma D^\alpha\big(\overline{\nabla}\cdot \frac{\overline{\nabla}\mathring \psi}{|\mathring N|}\big)
	+ \kappa D^\alpha(1-\overline{\Delta})^2\psi
	 + \kappa D^\alpha(1-\overline{\Delta})\partial_t\psi \big)(\overline{v}\cdot \overline{\nabla})  D^\alpha \mathring \psi
	-\int_\Sigma D^\alpha \mathring \psi \partial_3 q (\overline{v}\cdot \overline{\nabla})  D^\alpha \mathring \psi,
	\end{align*}
	where the last term can be directly controlled by $C(\mathring K_0)|\partial_3 q|_\infty|v|_\infty$. For the first term, we have 
	\begin{align*}
	\int_\Sigma  -\sigma D^\alpha\big(\overline{\nabla}\cdot \frac{\overline{\nabla}\mathring \psi}{|\mathring N|}\big)	D^\alpha \mathring \psi
	&\leq \sigma C(\mathring K_0,\kappa^{-1}),\\
	\int_\Sigma \kappa D^\alpha(1-\overline{\Delta})^2\psi(\overline{v}\cdot \overline{\nabla})  D^\alpha \mathring \psi 
	&= \int_\Sigma \kappa D^\alpha(1-\overline{\Delta})\psi(1-\overline{\Delta})\big((\overline{v}\cdot \overline{\nabla})  D^\alpha \mathring \psi\big)\\
	&\leq C(\mathring K_0)\int_\Sigma |\sqrt\kappa D^\alpha\psi |_2|v|_{W^{2,\infty}},\\
	\int_0^t\int_\Sigma\kappa D^\alpha(1-\overline{\Delta})\partial_t\psi (\overline{v}\cdot \overline{\nabla})  D^\alpha \mathring \psi
	&=\int_0^t\int_\Sigma\kappa D^\alpha\japbr\partial_t\psi \japbr\big((\overline{v}\cdot \overline{\nabla})  D^\alpha \mathring \psi\big) \\
	&\leq \delta\int_0^t\int_\Sigma|\sqrt\kappa\partial_tD^\alpha\psi|_1^2 +C(\mathring K_0)\int_0^t|v|^2_{W^{1,\infty}},
	\end{align*}
	where the control of the second term needs a minor modification when $D^\alpha =\partial_t^4$, since $\mathring K_0$ controls $\psi$ up to 7 derivatives but requires it to be attached with $4$ spatial derivative. In that case, we can exchange $\partial_t$ of $\mathring\psi$ with $\japbr$ of $\psi$ by integration by parts to close the estimate
	\begin{align*}
		&\int_0^t\int_\Sigma \kappa \partial_t^4(1-\overline{\Delta})^2\psi(\overline{v}\cdot \overline{\nabla})  \partial_t^4 \mathring \psi \\
		&= -\int_0^t\int_\Sigma \kappa \partial_t^5\japbr\psi(1-\overline{\Delta})\japbr\big((\overline{v}\cdot \overline{\nabla})  \partial_t^3 \mathring \psi \big)
		-\int_0^t\int_\Sigma \kappa \partial_t^4(1-\overline{\Delta})\psi(1-\overline{\Delta})\big((\partial_t\overline{v}\cdot \overline{\nabla})  \partial_t^3 \mathring \psi \big)\\
		&\quad +\int_\Sigma \kappa \partial_t^4\japbr\psi(1-\overline{\Delta})\big((\overline{v}\cdot \overline{\nabla})  \partial_t^3 \mathring \psi \big)\bigg|_0^t\\
		&\leq P(\mathring E_4^\kappa (0))+  \delta\int_0^t\int_\Sigma|\sqrt\kappa \partial_t^5\psi|_1^2 +C(\mathring K_0)\big(\int_0^t||v||_4^2+\int_0^t|\sqrt\kappa\partial_t^4\psi|_4+\int _\Sigma |\sqrt\kappa \partial_t^4 \psi|_1||v||_{3}\big).
	\end{align*}
	For the last term of (\ref{linpart2}), when $D^\alpha$ contains spatial derivatives, the analysis of $\int_\Sigma \mathring{\mathbf Q}\mathring{\mathcal{S}}_1$ and $\int _\Omega \mathring {\mathcal C}_i(v_i)\mathring{\mathbf Q}\partial_3\mathring\varphi$ is direct and identical to the one in section \ref{fullspatial}. We should only take a look at the case when $D^\alpha=\partial_t^4$, where cancellation structure occurs,
	\begin{align*}
	-\int_\Sigma \mathring{\mathbf Q}\mathring{\mathcal{S}}_1 
	&= 	\int_\Sigma \partial_t^4 q\mathring{\mathcal{S}}_1-\int_\Sigma \partial_t^4\mathring \varphi\partial_3^{\mathring \varphi} q\mathring{\mathcal{S}}_1,\\
	&= \int_\Sigma \partial_t^4 q \partial_t^4 \mathring \psi \partial _3 v\cdot \mathring N 
	-\int_\Sigma \partial_t^4 q[\partial_t^4,\overline{v}\cdot,\overline{\partial}\mathring\psi]-\int_\Sigma \partial_t^4\mathring \varphi\partial_3^{\mathring \varphi} q\mathring{\mathcal{S}}_1\\
	&=\int_\Sigma \partial_t^4 q\mathring \partial_t^4 \mathring \psi \partial _3 v\cdot \mathring N 
	-\int_\Sigma 4\partial_t^4 q \partial_t^3 \overline v\cdot  \partial_t \parbar\mathring\psi 
	-\int _\Sigma \partial_t^4q \sum _{k=1}^2 \begin{pmatrix} 4 \\ k \end{pmatrix} \partial_t^kv\cdot \partial_t^{4-k}\overline{\partial}\mathring \psi-\int_\Sigma \partial_t^4\mathring \varphi\partial_3^{\mathring \varphi} q\mathring{\mathcal{S}}_1.
	\end{align*}
	The second term is one that we had a cancellation, the rest are of lower order and can be directly controlled subject to the integration of $\partial_t$ by parts. The only trouble terms in $\int _\Omega \mathring {\mathcal C}_i(v_i)\mathring{\mathbf Q}\partial_3\mathring\varphi$ are the ones containing both $\partial_t^4 q$ and top order terms of $v$, for $\tau=1,2,$
	\begin{align*}
\mathring J_2^*&=4\int_\Omega \partial_t^4q \partial_t\big(\frac{\partial_\tau\mathring\varphi}{\partial_3 \mathring \varphi}\big)\partial_t^3\partial_3 v_\tau \partial_3\mathring \varphi =
4\int_\Omega \partial_t^4q (\partial_t\parbar\mathring\varphi\cdot \partial_t^3\partial_3 \overline{v}) 
-4\int_\Omega \partial_t^4q \partial_t \partial_3\mathring \varphi\partial_\tau\mathring\varphi\partial^{\mathring \varphi}_3\partial_t^3 v_\tau=:\mathring J^*_{21}+\mathring J^*_{22},\\
\mathring J_3^*&=-4\int_\Omega \partial_t^4q \partial_t\big(\frac{1}{\partial_3\mathring \varphi}\big)\partial_t^3 \partial_3 v_3 \partial_3\mathring \varphi =4\int_\Omega \partial_t^4q \partial_t\partial_3\mathring \varphi\partial^{\mathring \varphi}_3\partial_t^3  v_3  ,	
\end{align*}
where we can see the boundary integral produced by integrating $\partial_3$ in $\mathring J^*_{21}$ still cancels the trouble boundary integral $-\int_\Sigma 4\partial_t^4 q \partial_t^3 \overline v\cdot  \partial_t \parbar\mathring\psi$,
\begin{align*}
\mathring J_{21}^*	= -4\int_\Omega \partial_3\partial_t^4q (\partial_t\parbar\mathring\varphi\cdot \partial_t^3 \overline{v}) +4\int_\Sigma  \partial_t^4q (\partial_t\parbar\mathring\varphi\cdot \partial_t^3 \overline{v}),
\end{align*}
where the interior integral can be controlled by integrating $\partial_t$ by parts. Moreover, the top order term in $\mathring J_{22}^*+\mathring J^*_3$ still vanishes thanks to $\nabla^{\mathring \varphi}\cdot v =0$ and the remainder terms can be handled by parallel arguments in section \ref{fulltime}.
\\

\noindent\textbf{Summary}\\
We have derived the following uniform-in-$n$ estimate of the linearised system 
\begin{align}
\mathring{E}_4^\kappa(t)\leq \delta	\mathring{E}_4^\kappa(t)+C(\mathring K,\kappa^{-1})\bigg(\mathring{E}_4^\kappa(0)+\int_0^t\mathring{E}_4^\kappa\bigg),
\end{align}
where the $\delta	\mathring{E}_4^\kappa(t)$ can be absorbed to the energy on the left-hand side choosing small $\delta>0$. By Grownwall's inequality, we complete the proof of Prop \ref{higherorderlinearised}.

	\subsection{Picard iteration}
	So far, we have established the local well-posedness and the uniform-in-$n$ estimate for the linearised system (\ref{linearsys}) for each fixed $\kappa>0$. It suffices to prove $\{(q^{[n]}, v^{[n]}, F_k^{[n]},\psi^{[n]})\}$ has a strongly convergent subsequence, then the limit becomes the solution to the nonlinear approximate system (\ref{approxsys}). 
		
		For a function sequence $\{(f^{[n]})\}$, we define $[f]^{[n]} := f^{[n+1]}-f^{[n]}$, then we can write the linear system of $\{(q^{[n]}, v^{[n]}, F_k^{[n]},\psi^{[n]})\}$ as follows
		\begin{equation}\label{picardsys}
		\begin{cases}
			D_t^{\varphi^{[n]}} [v]^{[n]}- (\mathfrak F_k^{[n]}\cdot \nabla^{\varphi^{[n]}}) [F_k]^{[n]} +\nabla^{\varphi^{[n]}} [q]^{[n+1]}=-\mathring f_v^{[n]} \qquad & \text{in } \Omega,\\
			\nabla^{\varphi^{[n]}} \cdot [v]^{[n]}=-\mathring f^{[n]}_q & \text{in } \Omega,\\
			D_t^{\varphi^{[n]}} [F_j]^{[n]}-(\mathfrak F_j^{[n]}\cdot \nabla^{\varphi^{[n]}}) v^{[n]}=-\mathring f^{[n]}_F  & \text{in } \Omega,\\
			
			\partial_t [\psi]^{[n]}=[v]^{[n]}\cdot N^{[n]}+v^{[n]}\cdot [N]^{[n-1]} &\text{on } \Sigma,\\
			[q]^{[n+1]}=\sigma (\mathcal{H}(\psi^{[n]})-\mathcal{H}(\psi^{[n-1]})) + \kappa(1-\overline{\Delta})^2[\psi] ^{[n]}+ \kappa(1-\overline{\Delta})\partial_t[\psi]^{[n]} &\text{on } \Sigma,\\
			v_3^{[n+1]}=v_3^{[n]}=0&\text{on } \Sigma_b,\\
			([v]^{[n]},[F_k]^{[n]}, [q]^{[n]}, [\psi]^{[n]})|_{t=0} = (0,0,0,0).
		\end{cases}
	\end{equation}
	where the source terms are defined by 
	\begin{align}
		\mathring f_v^{[n]} &= [\overline{v}]^{[n-1]}\cdot \overline{\nabla}v^{[n]}+[ V_\mathbf N]^{[n-1]}\partial_3 v^{[n]} - ([\mathfrak{F}_k]^{[n-1]}\cdot \overline{\nabla})F_k^{[n]} - \mathfrak{F}_{k,\mathbf N}^{[n-1]}\partial_3 F_k^{[n]}+[\frac{\mathbf{N}}{\partial_3 \varphi}]^{[n-1]}\partial_3 q^{[n]},\\
		\mathring f_q^{[n]} &= [\frac{\mathbf{N}}{\partial_3 \varphi}]^{[n-1]}\cdot \partial_3 v^{[n]},\\
		\mathring f_F^{[n]}&= [\overline{v}]^{[n-1]}\cdot \overline{\nabla}F_k^{[n]}+[ V_\mathbf N]^{[n-1]}\partial_3 F_k^{[n]}- ([\mathfrak{F}_k]^{[n-1]}\cdot \overline{\nabla})v^{[n]} - \mathfrak{F}_{k,\mathbf N}^{[n-1]}\partial_3 v^{[n]},
	\end{align}
	with 
	\begin{align}
	V_\mathbf N^{[n]} =\frac{1}{\partial_3 \varphi^{[n]}}(v^{[n]}\cdot \mathbf N^{[n-1]}-\partial_t\varphi^{[n]}), \quad  \mathfrak{F}_{k,\mathbf N}^{[n]} := \frac{1}{\partial_3 \varphi^{[n]}}(F_k^{[n]}\cdot \mathbf N^{[n]}).
	\end{align}
	
	For $1\leq n$, we define the energy of the linear system (\ref{picardsys}) by 
	\begin{align}
	[E]_3^{[n]}(t): = 	\sum_{l=0}^3||\partial_t^l(  [v]^{[n]}, [F_k]^{[n]})||^2_{3-l}+ \sum_{l=0}^3|\sqrt\kappa \partial_t^l  [\psi]^{[n]}|^2_{5-l}+\int_0^t|\sqrt\kappa \partial_t^4 [\psi]^{[n]}|^2_1.
	\end{align}
	We shall prove the following estimate,
	\begin{prop}
		There exists $T_\kappa'>0 $ such that for $t\in [0,T_\kappa']$,
		\begin{align}
		[E]_3^{[n]}(t)\leq \frac{1}{4}\big(	[E]_3^{[n-1]}(t)+[E]_3^{[n-2]}(t)\big).
		\end{align}
	\end{prop}
	Note that $[E]_3^{[n-1]}(t)$ is produced by $[v]^{[n-1]},[F_k]^{[n-1]},[\psi]^{[n-1]}$ and $[E]_3^{[n-2]}$ is produced by $\mathbf N^{[n-2]}$ in the source terms.  The analysis is parallel to the estimate of $\mathring E_4(t)$ in section \ref{highlin} and we will only write the sketch of the proof.

	\subsubsection{Div-curl estimates}
  By lemma \ref{hodge}, for $l=0,1,2,$ we have 
 \begin{align}
		||\partial_t^l [v]^{[n]}||_{3-l}^2&\lesssim C(\mathring K_0)\big(||\nabla^{\varphi^{[n]}}\cdot \partial_t^l [v]^{[n]}	||_{2-l}+||\nabla^{\varphi^{[n]}}\times \partial_t^l [v]^{[n]}||_{2-l}^2+||\overline{\partial}^{3-l} \partial_t^l [v]^{[n]}||_0^2+||\partial_t^l[v]||_0^2\big),\\
		||\partial_t^l[F_j]^{[n]}||_{3-l}^2&\lesssim C(\mathring K_0)\big(||\nabla^{\varphi^{[n]}}\cdot \partial_t^l [F_j]^{[n]}||_{2-l}^2+||\nabla^{\varphi^{[n]}}\times \partial_t^l [F_j]^{[n]}||_{2-l}^2+||\overline{\partial}^{3-l} \partial_t^l [F_j]^{[n]}||_0^2+||\partial_t^l [F_j]^{[n]}||_0^2\big).
	\end{align}
	
	The $L^2$ estimate is straightforward and parallel to the arguments for the linear system (\ref{linearsys2}), we do not repeat it. The divergence part of $[v]^{[n]}$ can be directly controlled by $C(\mathring K_0)\int_0^t([E]_3^{[n]}+[E]_3^{[n-1]})$ invoking the second equation of (\ref{picardsys}). For the divergence part of $[F_j]^{[n]}$, we take $\nabla^{\varphi^{[n]}}\cdot$ to the evolution equation of $[F_j]^{[n]}$ to get
	\begin{align*}
	D_t^{\varphi^{[n]}}(\nabla^{\varphi^{[n]}}\cdot[F_j]^{[n]}) = 	(\partial^{\varphi^{[n]}}_i {\mathfrak F_j}^{[n]}\cdot \nabla^{\varphi^{[n]}}) [v]_i^{[n]}-(\mathfrak F_j^{[n]}\cdot \nabla^{\varphi^{[n]}})\mathring f_q^{[n]}+[D_t^{\varphi^{[n]}},\nabla^{\varphi^{[n]}}\cdot ][F_j]^{[n]},
	\end{align*}
	we observe that there are at most 1 derivative of $[v]^{[n]},[F_j]^{[n]},\mathring f_q^{[n]}$ and we can control $\mathring f_q^{[n]}$ up to 3 derivatives. The standard $H^2$ estimates and linearised transport theorem give us
	\begin{align}
	||\nabla^{\varphi^{[n]}}\cdot \partial_t^l [F_j]^{[n]}||^2_{2-l}
	&\leq C(\mathring K_0)(\nabla^{\varphi^{[n]}}\cdot \partial_t^l [F_j]^{[n]}(0)+\int_0^t[E]_3^{[n]}+[E]_3^{[n-1]}+[E]_3^{[n-2]})\\
	& = C(\mathring K_0)\int_0^t([E]_3^{[n]}+[E]_3^{[n-1]}+[E]_3^{[n-2]}).
	\end{align} 
	For the vorticity parts, taking $\nabla^{\varphi^{[n]}}\times$ to the evolution equations gives us
	\begin{align}
		&D_t^{\varphi^{[n]}}(\nabla^{\varphi^{[n]}}\times [v]^{[n]})\nonumber\\
		&=(\mathfrak F_k^{[n]}\cdot\nabla^{\varphi^{[n]}})(\nabla^{\varphi^{[n]}}\times [F_k]^{[n]})+(\nabla^{\varphi^{[n]}}\mathfrak F_{sk})\times (\partial_s^{\varphi^{[n]}}[F_k]^{[n]})+[D_t^{\varphi^{[n]}},\nabla^{\varphi^{[n]}}\times][v]^{[n]}-(\nabla^{\varphi^{[n]}}\times \mathring f_v^{[n]}),\\
		&D_t^{\varphi^{[n]}}(\nabla^{\varphi^{[n]}}\times [F_j]^{[n]})\nonumber\\
		&=(\mathfrak F_j^{[n]}\cdot\nabla^{\varphi^{[n]}})(\nabla^{\varphi^{[n]}}\times [v]^{[n]})+(\nabla^{\varphi^{[n]}}\mathfrak F_{sj})\times (\partial_s^{\varphi^{[n]}}[v]^{[n]})+[D_t^{\varphi^{[n]}},\nabla^{\varphi^{[n]}}\times][F_j]^{[n]}-(\nabla^{\varphi^{[n]}}\times \mathring f_F^{[n]}).
	\end{align}
	The top order terms $(\mathfrak F_k^{[n]}\cdot\nabla^{\varphi^{[n]}})(\nabla^{\varphi^{[n]}}\times [F_k]^{[n]})$ and $(\mathfrak F_j^{[n]}\cdot\nabla^{\varphi^{[n]}})(\nabla^{\varphi^{[n]}}\times [v]^{[n]})$ will be cancelled out in the energy estimate as we pointed out in section \ref{highlin}. Hence the standard $H^2$ estimates give us 
	\begin{align}
		||\nabla^{\varphi^{[n]}}\times \partial_t^l [v]^{[n]}||^2_{2-l}+||\nabla^{\varphi^{[n]}}\times \partial_t^l [F_k]^{[n]}||^2_{2-l}\leq C(\mathring K_0)\int_0^t([E]_3^{[n]}+[E]_3^{[n-1]}+[E]_3^{[n-2]}).
	\end{align}
	\subsubsection{Tangential estimates}
	Let $D^\alpha = \partial_t^{\alpha_0}\partial_1^{\alpha_1}\partial_2^{\alpha_2}$ with $|\alpha|\leq 3$. We define the AGUs of $[v],[F_j],[q]$ by 
	\begin{align*}
	[\mathbf V]^{[n]}:&= D^\alpha [v]^{[n]}-D^\alpha\varphi^{[n]}	\partial_3\varphi^{[n]}[v]^{[n]}, \quad [\mathbf F_j]^{[n]}:= D^\alpha [F_j]^{[n]}-D^\alpha\varphi^{[n]}	\partial_3\varphi^{[n]}[F_j]^{[n]},\\
	[\mathbf Q]^{[n]}&:=D^\alpha [q]^{[n]}-D^\alpha\varphi^{[n]}	\partial_3\varphi^{[n]}[q]^{[n]},
	\end{align*}

	\noindent the $D^\alpha$-differentiated (\ref{picardsys}) can be written in terms of the AGUs as follows

\begin{align}
\begin{cases}
	D_t^{\varphi^{[n]}}[\mathbf{V}]^{[n]}_i + \partial_i^{\varphi^{[n]}}  [{\mathbf{Q}}]^{[n]}= ({\mathfrak F}^{[n]}_k\cdot \nabla^{\varphi^{[n]}} ) [\mathbf{F}]^{[n]}_{ik}+\mathcal R_{v,i}^{[n]} \qquad & \text{in } \Omega,\\
	\partial_k^{\varphi^{[n]}}  [\mathbf{V}]^{[n]}_k=-\mathring{\mathcal{C}}_k([v]_k) & \text{in } \Omega,\\
	D_t^{\varphi^{[n]}}[\mathbf{F}]^{[n]}_{ij}= ({\mathfrak F}^{[n]}_k\cdot \nabla^{\varphi^{[n]}} ) [\mathbf{V}]^{[n]}_{ij}+\mathcal R_{F,ij}^{[n]}  & \text{in } \Omega, \\
	[\mathbf{Q}]^{[n]} = \sigma D^\alpha(\mathcal{H}(\psi^{[n]}-\mathcal{H}(\psi^{[n-1]}))
	+ \kappa D^\alpha(1-\overline{\Delta})^2[\psi]^{[n]} + \kappa D^\alpha(1-\overline{\Delta})\partial_t[\psi]^{[n]} 
	-D^\alpha \psi^{[n]} \partial_3 [q]^{[n]}   &\text{on }\Sigma , \\
	\partial_t D^\alpha [\psi]^{[n]} =[\mathbf{V}]^{[n]}\cdot N^{[n]}-(\overline{[v]}\cdot \overline{\nabla})  D^\alpha \psi^{[n]}+{\mathcal{S}}^{[n]}_b &\text{on }\Sigma. 
\end{cases} 
\end{align}
where 
\begin{align}
\mathcal{R}^{[n]}_{v,i}&=-D^\alpha \mathring f_v^{[n]}+\mathfrak F_{lk}^{[n]}\mathring{\mathcal{C}}_l([F]^{[n]}_{ik})+[D^\alpha, \mathfrak F^{[n]}_{lk}] \partial_l^{\varphi^{[n]}} [F]^{[n]}_{ik}-\mathring{\mathcal{D}}([v]_i)-\mathring{\mathcal{C}}_i([q]), \\
\mathcal R_{F,ij}^{[n]}&=-D^\alpha \mathring f_F^{[n]}+\mathfrak F^{[n]}_{kj}\mathring {\mathcal{C}}_k([v]_i)+[D^\alpha,\mathfrak F^{[n]}_{kj}]\partial_k^{\varphi^{[n]}} [v]_i-\mathring{\mathcal{D}}([F]_{ij}),\\
\mathcal{S}_b^{[n]}&=D^\alpha \psi^{[n]} \partial _3 [v]\cdot  N^{[n]} -[D^\alpha,\overline{[v]}^{[n]}\cdot,\overline{\partial}\psi^{[n]}]+D^\alpha(v^{[n]}\cdot [N]^{[n-1]}).
\end{align}
Note that the additional remainder terms can be directly controlled by $C(\mathring K_0)\int_0^t( [E]_3^{[n-1]}+[E]_3^{[n-2]})$, hence, we can apply analysis similar to the ones in (\ref{lintan}) to obtain
\begin{align*}
	&\frac{1}{2}\frac{d}{dt}\int_\Omega \bigg(|[\mathbf V]_i^{[n]}|^2+|[\mathbf F_k]^{[n]}|^2\partial_3\varphi^{[n]}\bigg)\\
	&\leq C(\mathring K_0 )\int_0^t\bigg([E]_3^{[n]}+[E]_3^{[n-1]}+[E]_3^{[n-2]}\bigg)-\int_\Sigma [\mathbf Q]^{[n]}([\mathbf V]^{[n]}\cdot N^{[n]})+\int _\Omega D^\alpha [q]^{[n]}\mathring {\mathcal {C}}_i([v]_i)\partial_3 \varphi^{[n]}.
\end{align*}
We expand the boundary integral
\begin{align*}
-\int_\Sigma [\mathbf Q]^{[n]}([\mathbf V]^{[n]}\cdot N^{[n]}) = 	
-\int_\Sigma [\mathbf Q]^{[n]}\partial_t D^\alpha [\psi]^{[n]}
-\int_\Sigma [\mathbf Q]^{[n]}(\overline{[v]}^{[n]}\cdot \overline{\nabla})  D^\alpha \psi^{[n]}
+\int_\Sigma [\mathbf Q]^{[n]}{\mathcal{S}}^{[n]}_b ,
\end{align*}
where the first boundary integral still contributes to the boundary regularity
\begin{align*}
	&-\int_0^t\int_\Sigma [\mathbf Q]^{[n]}\partial_t D^\alpha [\psi]^{[n]}+\int_\Sigma|\sqrt\kappa (1-\overline{\Delta})D^\alpha[\psi]^{[n]}|^2\bigg|_0^t+\int_0^t\int_\Sigma |\sqrt{\kappa}\japbr\partial_tD^\alpha [\psi]^{[n]}|^2\\
	&\leq C(\mathring K_0)\int_0^t \bigg([E]_3^{[n]}+[E]_3^{[n-1]}\bigg).
\end{align*}
The second term can be directly bounded by the $\kappa-$weighted regularity invoking the boundary condition for $[\mathbf Q]^{[n]}$ and integrating $(1-\overline{\nabla}),\japbr$ by parts. For $\int_\Sigma [\mathbf Q]^{[n]}{\mathcal{S}}^{[n]}_b+\int _\Omega D^\alpha [q]^{[n]}\mathring {\mathcal {C}}_i([v]_i)\partial_3 \varphi^{[n]}$, we can direct control it when $D^\alpha$ contains spatial derivatives. When $D^\alpha =\partial_t ^4$, we can still find the same cancellation structure that the troubling term $-4\int_\Sigma \partial_t ^3 [q]^{[n]}\partial_t^2\overline{[v]}^{[n]}\partial_t\parbar\psi^{[n]} $ from $\int_\Sigma [\mathbf Q]^{[n]}{\mathcal{S}}^{[n]}_b$ can still be cancelled by the boundary integral produced by integrating $\partial_3$ by parts in $4\int_\Omega\partial_t^3[q](\partial_t\parbar\varphi^{[n]}\cdot\partial_t^2\partial_3\overline{[v]})\partial_3 \varphi^{[n]}$ from $\int _\Omega D^\alpha [q]^{[n]}\mathring {\mathcal {C}}_i([v]_i)\partial_3 \varphi^{[n]}.$

Summarising the above estimates, for small $\delta>0$, we obtain the energy estimate 
\begin{align}
[E]_3^{[n]}(t)\leq \delta	[E]_3^{[n]}(t)+C(\mathring K_0,\kappa^{-1})\int_0^t \bigg([E]_3^{[n]}+[E]_3^{[n-1]}+[E]_3^{[n-2]}\bigg).
\end{align}
Thus there exists $T_\kappa '>0$ independent of $n$ such that
\begin{align}
	\sup_{0\leq t\leq T_\kappa'}[E]_3^{[n]}(t)\leq \frac{1}{4}\bigg(\sup_{0\leq t\leq T_\kappa'}[E]_3^{[n-1]}(t)+\sup_{0\leq t\leq T_\kappa'}[E]_3^{[n-2]}(t)\bigg),
\end{align}
and thus we know by induction that
\begin{align}
\sup_{0\leq t\leq T_\kappa'}[E]_3^{[n]}(t)\leq C(\mathring K_0,\kappa^{-1})\frac{1}{2^{n-1}}\to 0 \text{ as } n\to \infty	.
\end{align}
Hence, for fixed $\kappa >0 $, the sequence of approximate solutions $\{(v^{[n]},F_k^{[n]},\mathfrak F_k^{[n]},q^{[n]},\psi^{[n]})\}$ has a strongly convergent subsequence and the limit $\{(v^{[\infty]},F_k^{[\infty]},\mathfrak  F_k^{[\infty]},q^{[\infty]},\psi^{[\infty]})\}$ is a solution of (\ref{limitsys}). We can recover the solution $\{(v^{[\infty]},F_k^{[\infty]},q^{[\infty]},\psi^{[\infty]})\}$ to our approximate system (\ref{approxsys}) by remark \ref{recover}.

\pagebreak

\bibliographystyle{plain}
\bibliography{LWP.bib}

\end{document}